\documentclass[nonblindrev]{informs3_noinforms}              

\usepackage[numbers]{natbib}

\usepackage{enumerate}
\newtheorem{observation}{Observation}
\newtheorem{Question}{Question}

\newcommand{\pr}{\mathbb{P}}
\newcommand{\E}{\mathbb{E}}

\newcommand{\reals}{{\mathbb R}}
\newcommand{\ignore}[1]{\relax}

\def\bbr{{\Bbb{R}}} 

\renewcommand{\theequation}{\thesection.\arabic{equation}}

\usepackage{natbib}
 \NatBibNumeric
 \def\bibfont{\small}%
 \bibpunct[, ]{[}{]}{,}{n}{}{,}%

\usepackage[colorlinks=true,breaklinks=true,bookmarks=true,urlcolor=blue,
     citecolor=blue,linkcolor=blue,bookmarksopen=false,draft=false]{hyperref}

\def\EMAIL#1{\href{mailto:#1}{#1}}
\def\URL#1{\href{#1}{#1}}         

\TheoremsNumberedThrough     

\EquationsNumberedThrough    

\begin{document}



\RUNTITLE{On the steady-state probability of delay and large negative deviations in the Halfin-Whitt regime}

\TITLE{On the steady-state probability of delay and large negative deviations for the $GI/GI/n$ queue in the Halfin-Whitt regime}

\ARTICLEAUTHORS{%
\AUTHOR{\bf David A. Goldberg}
\AFF{Georgia Institute of Technology, \EMAIL{dgoldberg9@isye.gatech.edu}, \URL{http://www2.isye.gatech.edu/~dgoldberg9/}}
}
\RUNAUTHOR{{David A. Goldberg}}


\ABSTRACT{%
We consider the FCFS $GI/GI/n$ queue in the Halfin-Whitt heavy traffic regime, and prove bounds for the steady-state probability of delay (s.s.p.d.) for generally distributed processing times.  We prove that there exist $\epsilon_1, \epsilon_2 > 0$, depending on the first three moments of the inter-arrival and processing time distributions, such that the s.s.p.d. is bounded from above by $\exp\big(-\epsilon_1 B^2\big)$ as the associated excess parameter $B \rightarrow \infty$; and by $1 - \epsilon_2 B$ as $B \rightarrow 0$.  We also prove that the tail of the steady-state number of idle servers has a Gaussian decay.  We provide explicit bounds in all cases, in terms of the first three moments of the inter-arrival and service distributions, and use known results to show that our bounds correctly capture various qualitative scalings.
\\\indent Our main proof technique is the derivation of new stochastic comparison bounds for the FCFS $GI/GI/n$ queue, which are of a structural nature, hold for all $n$ and times $t$, and significantly generalize the work of \citet{GG.10c} (e.g. by providing bounds for the queue length to exceed any given level, as opposed to any given level strictly greater than the number of servers as acheived in \citet{GG.10c}).  Our results do not follow from simple comparison arguments to e.g. infinite-server systems or loss models, which would in all cases provide bounds in the opposite direction. 

}%


\KEYWORDS{many-server queues, Halfin-Whitt regime, probability of delay, stochastic comparison, weak convergence, large deviations, Gaussian process, Slepian's inequality, renewal process}

\MSCCLASS{60K25}

\maketitle

%


\section{Introduction.}\label{Section:Introduction}
Parallel server queueing systems can operate in a variety of regimes that balance between efficiency and quality of offered service.  This is captured by the so-called Halfin-Whitt (H-W) heavy traffic regime.  Although studied originally by \citet{E.48} and \citet{J.74}, the regime was formally introduced by \citet{HW.81}, who studied the $GI/M/n$ system (for large $n$) when the traffic intensity $\rho$ scales like $1 - Bn^{-\frac{1}{2}}$ for some strictly positive excess parameter $B$.  The authors prove that for this sequence of $GI/M/n$ queueing models (indexed by $n$), the steady-state probability that an arriving job has to wait for service (i.e. steady-state probability of delay) converges (as $n \rightarrow \infty$) to a function of $B$ (independent of $n$), which they explicitly compute.  This limiting probability converges to $0$ as $B \rightarrow \infty$ (low-utilization regime), converges to $1$ as $B \rightarrow 0$ (high-utilization regime), and decreases monotonically from $1$ to $0$ as $B$ increases from $0$ to $\infty$, thus nicely quantifying the trade-off between server utilization and quality of service (as measured by s.s.p.d.).  Analogous explicit formulas have been found for the case in which processing times are a mixture of an exponential distribution and a point mass at 0 (i.e. $H^*_2$ distributed)  by \citet{Whitt.05d}; and for the case of deterministic processing times by \citet{JMM.04} (and \citet{hurvich2014seriesb} in the presence of heavy-tailed inter-arrivals).  For the case of exponentially distributed processing times, these results have also been extended to allow for abandonments, see \citet{GMR.02d} and \citet{ZM.05d}.
\\\indent However, much less is known for more general processing time distributions.  This is particularly unfortunate, as it is believed that in many of the major application domains of the H-W regime, processing times have more general distributions (e.g. log-normal), see \citet{Mandelbaum.05}.  In this more general setting, the main known results with regards to (w.r.t.) the s.s.p.d. can be summarized as follows.  For processing times with finite support, \citet{GM.08} give an implicit description (in terms of a certain Markov chain) of the limiting s.s.p.d in the H-W regime (also proving that this limit exists), and show that this probability lies strictly in $(0,1)$.  \citet{Whitt.04d} gives several heuristic approximations for this limiting s.s.p.d. for generally distributed processing times, which he verifies numerically.  \citet{dai2013validity} prove the validity of a certain diffusion approximation for the steady-state number in system for the case of phase-type processing times with abandonments, and the work of \citet{he2011many} provides numerical algorithms for evaluating the associated probabilities.  
\\\indent A similar situation exists for questions related to the large deviations (l.d.) properties of the steady-state number of \emph{idle} servers in the H-W regime.  Indeed, the settings in which this limiting l.d. behavior is precisely understood coincide exactly with the settings in which the limiting s.s.p.d. can be computed explicitly, as described above.  Furthermore, very little is known about even the qualitative behavior / crude asymptotic scaling of either quantity.  However, even from the special cases of $GI/H^*_2/n$ and $GI/D/n$ queues for which the relevant limits can be computed in closed form, some interesting patterns emerge.  We now briefly describe these patterns heuristically before introducing all relevant model details, as a more complete discussion follows in Section\ \ref{4comparesec}.  Let $p_B$ denote the limiting s.s.p.d. when the excess parameter is $B$, and $F_{\textrm{idle},B}(x)$ denote the limiting steady-state probability of seeing more than $x n^{\frac{1}{2}}$ idle servers in the $n$th system under the H-W scaling, where the particular inter-arrival and processing time distributions are to be inferred from context (and for the moment restricting to inter-arrival and processing time distributions for which these limits are known to exist).  Note that as $B \rightarrow 0$, one would expect $p_B$ to converge to 1, since as $B \rightarrow 0$ the system becomes more overloaded.  Similarly, as $B \rightarrow \infty$, one would expect $p_B$ to converge to 0.  We have already asserted that \citet{HW.81} proved this for the case of exponentially distributed processing times.  However, one can ask \emph{how quickly} these quantities approach 1 as $B \rightarrow 0$, and $0$ as $B \rightarrow \infty$.  Such qualitative information would provide valuable insight into the s.s.p.d. in the H-W regime, beyond the initial insight of Halfin and Whitt that the correct limits should be 1 and 0.  We believe such insight to be especially relevant for queues in the H-W regime, since one of the primary motivating features of the H-W regime was that the s.s.p.d. should have a non-trivial limit as $n \rightarrow \infty$.
\\\indent Then the work of \citet{HW.81}, \citet{Whitt.05d}, and \citet{JMM.04} imply the following.  
\begin{observation}\label{obsintro}
\ \ \ When the processing times are exponentially distributed, $H^*_2$, or deterministic, and the inter-arrival distribution has finite third moment, there exist strictly positive finite constants $\epsilon_1,\epsilon_2,\epsilon_3$, depending only on the inter-arrival and processing time distributions (but not $B$ or $n$), such that (s.t.):
\begin{enumerate}[(i)]
\item $\lim_{B \rightarrow \infty} B^{-2} \log(p_B) = - \epsilon_1$, \label{introa}
\item $\lim_{B \rightarrow 0} B^{-1} (1 - p_B)  = \epsilon_2$, \label{introb}
\item $\lim_{x \rightarrow \infty} x^{-2} \log\big(F_{\textrm{idle},B}(x)\big) = - \epsilon_3$. \label{introc}
\end{enumerate}
\ \\Intuitively, this implies that $p_B$ behaves (roughly) like $\exp(- \epsilon_1 B^2)$ as $B \uparrow \infty$, and $1 - \epsilon_2 B$ as $B \downarrow 0$; and that $F_{\textrm{idle},B}(x)$ behaves like $\exp(- \epsilon_3 x^2)$ as $x \uparrow \infty$.
\end{observation}
As previously mentioned, it is believed that for the relevant systems which arise in practice, the service time distributions are in fact closer to a log-normal distribution, see \citet{Mandelbaum.05}.  In this paper, we thus set out to answer the following question.
\begin{Question}\label{CentralQ}
\ \ \ To what extent do the scalings suggested by Observation\ \ref{obsintro} hold in general?
\end{Question}
We note that the scaling of (\ref{introa}) and (\ref{introb}) represent ``limits within limits", where letting $B \rightarrow 0$ is moving the H-W scaling ``closer" to the classical heavy-traffic scaling in which $\rho$ scales like $1 - \frac{c}{n}$ for some $c > 0$ as $n \rightarrow \infty$, while letting $B \rightarrow \infty$ is moving the H-W scaling closer to the classical heavy-traffic scaling in which $\rho$ is bounded away from 1 as $n \rightarrow \infty$, and we refer the reader to \citet{HW.81} for a more detailed discussion.  We note that the regime in which $B \rightarrow 0$ has been explicitly discussed in the context of $GI/D/n$ queues by \citet{JMM.04}, \citet{janssen2007cumulants}, and \citet{hurvich2014seriesb}.  This setting is also related to the question of whether various quantities are differentiable w.r.t. the drift parameter, see \citet{dieker2014sensitivity} and \citet{lipshutz2016directional}.  We further note that if one wishes to bound the s.s.p.d. from below, or the number of idle servers from above (as done by \citet{GS.12d} Theorem 4.(i), and Lemmas 5 and 10), one can derive non-trivial bounds in the H-W regime by comparing the true system to an appropriate infinite-server queue, even for the case of generally distributed processing times, and we include a more complete discussion in Section\ \ref{4comparesec}.  However, it seems that no similarly straightforward comparison yields meaningful bounds in the other direction, i.e. upper bounds on the s.s.p.d. and lower bounds on the number of idle servers, in the H.W. regime.  In this paper, we will focus on exactly these settings where such arguments seem to break down.
\\\indent We now briefly review the more general question of when the limits $p_B$ and $F_{\textrm{idle},B}(x)$ are known to exist (for any fixed $B,x$), i.e. when the appropriately scaled sequence of queues under the H-W scaling has a weak limit, irregardless of the more subtle question of whether those limits themselves satisfy (\ref{introa}) - (\ref{introc}) as one varies $B$ and $x$.  As previously discussed, the relevant weak limits are proven to exist, and actually computed explicitly, for the case of exponentially distributed processing times by \citet{HW.81}, $H^*_2$ processing times by \citet{Whitt.05d}, and deterministic processing times by \citet{JMM.04}.  Similar, albeit often less explicit, weak convergence results under the H-W scaling were subsequently obtained for more general multi-server systems by \citet{R.00}, \citet{MM.08}, \citet{GM.08}, \citet{GS.12d}, \citet{kaspi2013spde}, \citet{R.09}, and \citet{PR.10}.  Many of these results have also been extended to the setting of abandonments by \citet{MM.12d}, \citet{DHT.10}, \citet{RT.12b}, \citet{huang2012diffusion}, \citet{dai2013validity}, \citet{kang2012asymptotic}, \citet{weerasinghe2013diffusion}, \citet{dai2014validity}, and we refer the interested reader to the surveys of \citet{DH.12g} and \citet{ward2012asymptotic}.  Using Stein's method, several recent papers have also provided explicit bounds on the associated rate of weak convergence, see \citet{gurvich2014diffusion}, \citet{braverman2015stein}, and \citet{braverman2016high}.  
\\\indent Unfortunately, as the theory of weak convergence generally relies heavily on the assumption of compact time intervals, the most general of these results hold only in the transient regime.  Indeed, the only settings in which the relevant sequence of normalized steady-state queue-length distributions have also been shown to have a weak limit are subsumed by the settings in which the processing time distribution either has finite support, or is of phase-type.  Although any distribution can be approximated to within any accuracy by either of these families, it is not known how the quality of such an approximation behaves under the H-W scaling.  We note that recent work on the so-called s.p.d.e. approach has made considerable progress towards proving existence of the relevant weak limits in much greater generality, see \citet{aghajani2015ergodicity} and \citet{aghajani2015ergodicityb}.  However, it is important to point out that the resulting limiting processes seem to have a complicated behavior (see \citet{PR.10}, \citet{dieker2013positive}, \citet{aghajani2015ergodicity}), and it is not clear how to implement a direct analysis of these processes to tackle Question\ \ref{CentralQ}.  
\begin{observation}\label{importantobs}
\ \ \ In spite of the extensive body of work on queues in the H-W regime, the answer to Question\ \ref{CentralQ} is unknown beyond the cases of exponentially distributed, $H^*_2$, and deterministic processing times.
\end{observation}
Recently, \citet{GG.10c} developed a novel stochastic comparison framework to address some of the issues faced by previous techniques, namely the difficulty of proving results in the steady-state, and the non-explicit descriptions of relevant limits.  In that work, the authors analyzed the steady-state $GI/GI/n$ queue by comparing to a ``modified system" in which all servers are kept busy at all times by adding ``artificial arrivals" to the system whenever a server would otherwise have gone idle.  They showed that the steady-state distribution of this modified system has a very simple representation as the supremum of a certain one-dimensional random walk.  Using these techniques, they proved tightness of the relevant sequence of normalized steady-state measures, and computed the limiting large deviation behavior for the number of jobs waiting in queue, under quite general assumptions.  Unfortunately, the techniques developed by \citet{GG.10c} cannot be used to provide upper bounds on the s.s.p.d., or lower bounds on the number of idle servers, as the aforementioned modified system always has all servers busy.  We note that although stochastic comparison techniques were also used in \citet{GS.12d}, in conjunction with Lyapunov function arguments, to analyze certain quantities in the H-W regime, it seems their approach similarly does not work for the quantities considered in Question\ \ref{CentralQ}.
\subsection{Our contribution.}
In this paper, we make considerable progress towards resolving Question\ \ref{CentralQ}.  We prove that under quite general assumptions on the inter-arrival and processing time distributions, e.g. finite third moment, there exists some strictly positive $\epsilon_1,\epsilon_2,\epsilon_3$ (depending on the first three moments of the inter-arrival and processing time distributions) s.t. in the H-W regime, the s.s.p.d. is bounded from above by $\exp\big(-\epsilon_1 B^2\big)$ as the associated excess parameter $B \rightarrow \infty$; and by $1 - \epsilon_2 B$ as $B \rightarrow 0$.  We also prove that the probability of there being more than $x n^{\frac{1}{2}}$ idle servers (in steady-state, for large $n$) is bounded from below by $\exp\big( - \epsilon_3 x^2 \big)$ as $x \rightarrow \infty$.  Combined with known results for the the $M/D/n, M/H^*_2/n$, and $M/GI/\infty$ queues, our results show that the qualitative scaling suggested by (\ref{introa}) - (\ref{introc}) indeed holds in considerable generality, yielding a partial positive answer to Question\ \ref{CentralQ}.  Furthermore, we provide explicit estimates for $\epsilon_1,\epsilon_2,\epsilon_3$, and use known results to prove that our lower (upper) bounds for $\epsilon_1$ ($\epsilon_3$) are close to known upper (lower) bounds for these quantities.  We note that although our proofs do not demonstrate the existence of the relevant limits, they do show that all relevant $\liminf, \limsup$, and associated subsequential limits scale in the aforementioned manner.
\\\indent Our main proof technique is the derivation of new and simple bounds for the FCFS $GI/GI/n$ queue, which significantly extend the stochastic comparison framework developed by \citet{GG.10c}.  We consider a different modified system, in which we not only keep certain servers busy by adding ``artificial arrivals", but also allow servers to ``break down" at a time of our choosing, leaving a number of ``working servers" which is strictly less than $n$.  By then keeping only those ``working servers" busy at all times, it becomes possible for the number of jobs in this second modified system to go below $n$, yielding meaningful bounds on the s.s.p.d. and number of idle servers in the H-W regime.  Our techniques allow us to analyze many properties of the $GI/GI/n$ queue in the H-W regime without having to consider the complicated exact dynamics of the $GI/GI/n$ queue.  Our results can also be viewed as a step towards developing a calculus of stochastic-comparison type bounds for parallel server queues, in which one derives bounds by composing structural modifications to a parallel-server queue (e.g. adding jobs, removing servers, adding servers, etc.) over time.  Although stochastic comparison techniques have been widely used to study queueing systems, see \citet{Stoyan.83}, we believe the bounds developed in this paper to be novel, and particularly suited to studying queues in the H-W regime.  Another appealing feature of our results is that, as in the work of \citet{GG.10c}, our bounds are of a structural nature, hold for all $n$ and all times $t \geq 0$, and have intuitive closed-form representations as the suprema of certain natural processes which converge weakly to Gaussian processes.
\subsection{Outline of paper.}
The rest of the paper proceeds as follows.  In Section\ \ref{4mainsec}, we present our main results.  In Section\ \ref{4uppersec}, we establish our general-purpose upper bounds for the queue length in a properly initialized FCFS $GI/GI/n$ queue.  In Section\ \ref{4upperasymptoticsec}, we prove an asymptotic version of our upper bound in the H-W regime.  In Section\ \ref{Binfinity}, we prove our bounds on the s.s.p.d. as $B \rightarrow \infty$.  In Section\ \ref{Bzero}, we prove our bounds on the  s.s.p.d. as $B \rightarrow 0$.  In Section\ \ref{negx}, we prove our bounds on the large deviations behavior of the steady-state number of idle servers.  In Section\ \ref{4comparesec}, we compare to previous results from the literature, which show that our bounds are tight in an appropriate sense.  In Section\ \ref{4concsec} we summarize our main results and comment on directions for future research.  We include a technical appendix in Section\ \ref{4appsec}.
\section{Main results.}\label{4mainsec}
We consider the First-Come-First-Serve (FCFS) $GI/GI/n$ queueing model, in which inter-arrival times are independent and identically distributed (i.i.d.) r.v.s, and processing times are i.i.d. r.v.s.
\\\indent Let $A$ and $S$ denote some fixed r.v.s with non-negative support s.t. $\E[A] = \mu_A^{-1} < \infty, \E[S] = \mu_S^{-1} < \infty$, and $\pr(A = 0) = \pr(S = 0) = 0$.  Let $\sigma_A$ and $\sigma_S$ denote the standard deviations of $A$ and $S$, respectively.  Also, let $c_A$ and $c_S$ denote the coefficient of variation of $A$ and $S$, respectively, i.e. $c_A = \mu \sigma_A$ and $c_S = \mu \sigma_S$.
\\\indent For excess parameter $B > 0$ and number of servers $n \geq 1$, let $\lambda_{n,B} \stackrel{\Delta}{=} n - B n^{\frac{1}{2}}$.
For $n$ sufficiently large to ensure $\lambda_{n,B} > 0$ (which is assumed throughout), let $Q^n_B(t)$ denote the number in system (number in service $+$ number waiting in queue) at time $t$ in the FCFS $GI/GI/n$ queue with inter-arrival times drawn i.i.d. distributed as $A\lambda_{n,B}^{-1}$ and processing times drawn i.i.d. distributed as $S$ (initial conditions will be specified later), independently from the arrival process.  Note that this scaling is analogous to that studied by \citet{HW.81}, as the traffic intensity in the $n$th system is $1 - B n^{-\frac{1}{2}}$ in both settings (supposing $\mu_A = \mu_S$).  All processes should be assumed right-continuous with left limits (r.c.l.l.) unless stated otherwise.  All empty summations should be evaluated as zero, all empty products should be evaluated as one, and all logarithms should be taken base $e$ unless specified otherwise.
\subsection{Main results.}
Our main results will require two additional sets of assumptions on $A$ and $S$.  The first set of assumptions, which we call the H-W assumptions, ensures that $\lbrace Q^n_B(t), n \geq 1 \rbrace$ is in the H-W scaling regime as $n \rightarrow \infty$.  We say that $A$ and $S$ satisfy the H-W assumptions iff $\mu_A = \mu_S$, in which case we denote this common rate by $\mu$.
The second set of assumptions, which we call the $T_0$ assumptions, is a set of additional technical conditions we require for our main results.
\begin{enumerate}[(i)]
\item There exists $\epsilon > 0$ s.t. $\E[A^{2+\epsilon}], \E[S^{2+\epsilon}] < \infty$.\label{t1}
\item $\sigma^2_A + \sigma^2_S > 0$.  Namely either $A$ or $S$ is a non-degenerate r.v.\label{t2}
\item $\limsup_{t \downarrow 0} t^{-1} \pr ( S \leq t ) < \infty$.\label{t3}
\item For all $B > 0$ and sufficiently large $n$, and all initial conditions, $Q^n_B(t)$ converges weakly to a stationary measure $Q^n_B(\infty)$ as $t \rightarrow \infty$, independent of initial conditions.\label{t4}
\end{enumerate}
We refer the interested reader to \citet{GG.10c} for a discussion of the restrictiveness and necessity of these assumptions, as that work uses the same set of assumptions.
\\\indent We now state our main results.  We begin by stating our bound on the s.s.p.d. as $B \rightarrow \infty$.
\begin{theorem}\label{largebtheorem} 
For any fixed $A$ and $S$ which satisfy the $H-W$ and $T_0$ assumptions,
$$\limsup_{B \rightarrow \infty} B^{-2} \log \bigg( \limsup_{n \rightarrow \infty} \pr \big( Q^n_B(\infty) \geq n \big) \bigg) \leq 
- \frac{1}{16} (2 c^2_S + c^2_A + 1)^{-1}.$$
\end{theorem}
Intuitively, Theorem\ \ref{largebtheorem} implies that in the H-W regime, the s.s.p.d. is bounded from above by $\exp\big(- \frac{1}{16} (2 c^2_S + c^2_A + 1)^{-1} B^2 \big)$ as $B \rightarrow \infty$.  
As we will see in Section\ \ref{4comparesec}, known results imply that the bounds of Theorem\ \ref{largebtheorem} are close to the best possible.
We now give our bounds on the s.s.p.d. as $B \rightarrow 0$.  Let $\alpha_S \stackrel{\Delta}{=} \mu^3 \bigg( \E[S^2] + 2 \E[S^3] + \frac{3}{8} \mu \big(\E[S^2]\big)^2 \bigg)$, and 
$$\delta_{A,S} \stackrel{\Delta}{=} \bigg(200 + \alpha_S + \mu^2 + \mu^{-2} + \sigma_S^2 + \sigma_S^{-2} + \sigma_A^2 + \sigma_A^{-2} + c_A^2 + c_A^{-2} + c_S^2 + c_S^{-2}\bigg)^{80}.$$
\begin{theorem}\label{smallbtheorem} 
For any fixed $A$ and $S$ which satisfy the $H-W$ and $T_0$ assumptions, s.t. in addition $\E[S^3] < \infty$, $\sigma^2_A > 0$, and $\sigma^2_S > 0$, 
\begin{equation}\label{smallbbb}
\liminf_{B \rightarrow 0} B^{-1} \liminf_{n \rightarrow \infty} \pr \big( Q^n_B(\infty) < n \big) \geq
\exp\big( - 10^{27} \delta_{A,S} \big).
\end{equation}
\end{theorem}
Intuitively, Theorem\ \ref{smallbtheorem} implies that the s.s.p.d. is bounded from above by $1 - \exp\big( - 10^{27} \delta_{A,S} \big) B$ as $B \rightarrow 0$.
Equivalently, the steady-state probability that a job \emph{does not} have to wait for service, i.e. no delay, is bounded from below by $\exp\big( - 10^{27} \delta_{A,S} \big) B$ as $B \rightarrow 0$.
We note that in contrast to our other results, Theorem\ \ref{smallbtheorem} is not a statement about the logarithm of a given probability, but of the probability itself.  For this reason, all constants involved in various estimates appearing in the relevant proofs must be incorporated into the final result.  Combined with the fact that when $B$ is small extra care must be taken to ensure that other terms do not dominate the relevant probabilities, our results for the setting in which $B \downarrow 0$ have an admittedly massive prefactor.  For this reason, we view the main insight of Theorem\ \ref{smallbtheorem} to be that the given probability is at most $1 - \epsilon_2 B$ for some strictly positive $\epsilon_2$ as $B \downarrow 0$, and that an explicit bound on $\epsilon_2$ depending only on the first three moments is acheivable.  Furthermore, as we will see in Section\ \ref{4comparesec}, known results imply that there exists $\epsilon_2' > 0$ s.t. the given probability is at least $1 - \epsilon_2' B$ as $B \downarrow 0$, at least for the setting of Markovian and deterministic processing times.
\\\indent We now state our bounds on the large deviations behavior for the number of idle servers.
\begin{theorem}\label{negxtheorem}
For any fixed $A$ and $S$ which satisfy the $H-W$ and $T_0$ assumptions, s.t. in addition $\sigma^2_A > 0$, and any fixed $B > 0$,
$$
\liminf_{x \rightarrow \infty} x^{-2} \log \Bigg( \liminf_{n \rightarrow \infty} \pr\bigg( \big( Q^n(\infty) - n \big)n^{-\frac{1}{2}} \leq - x \bigg) \Bigg) \geq - 2 \E[S^2] \sigma_A^{-2}.
$$
\end{theorem}
Intuitively, Theorem\ \ref{negxtheorem} implies that the tail of the limiting steady-state number of idle servers is bounded from below by $\exp\big( - 2 \E[S^2] \sigma_A^{-2} x^2 \big)$ as $x \rightarrow \infty$.  As we will see in Section\ \ref{4comparesec}, known results imply that the bounds of Theorem\ \ref{negxtheorem} are close to the best possible in many settings.
\\\indent In addition to our main results Theorems\ \ref{largebtheorem} - \ref{negxtheorem}, we note that our stochastic comparison approach implies novel bounds for the probabilities associated with a $GI/GI/n$ queue with any given inter-arrival and service distribution, irregardless of whether the system is in the H-W regime (see e.g. Theorem\ \ref{4ubound1}).  Furthermore, our approach also implies bounds for the s.s.p.d., and more generally any probability of interest, for the $GI/GI/n$ queue in the H-W regime for any given $B > 0$ (see e.g. Theorem\ \ref{mainasymptotic}).
Here we chose to focus our main results on Theorems\ \ref{largebtheorem}\ -\ \ref{negxtheorem} instead of Theorems\ \ref{4ubound1} - \ref{mainasymptotic} for clarity of exposition, as the latter involve a somewhat complicated 2-dimensional optimization over a nested supremum.  We do note that a careful analysis of Gaussian processes, similar to that which we use in the proofs of Theorems\ \ref{largebtheorem} - \ref{negxtheorem}, could likely be combined with Theorem\ \ref{mainasymptotic} to give relatively simple non-trivial bounds for any probability of interest associated with the $GI/GI/n$ queue in the H-W regime (for any given $B > 0$), although we do not pursue such an analysis here.
\\\indent We further note that our results cannot be derived using straightforward comparison bounds, to e.g. an infinite-server or loss system, as in all cases our inequalities point in the other direction.  Also, as in \citet{GG.10c}, our results translate into bounds for any weak limits of the associated sequences of r.v.s.
\section{Upper bound.}\label{4uppersec}
In this section, we prove general upper bounds for the FCFS $GI/GI/n$ queue, when properly initialized.  The bounds are valid for all finite $n$, and work in both the transient and steady-state (when it exists) regimes.  Although we will later customize these bounds to the H-W regime to prove our main results, we note that the bounds are in no way limited to that regime. 
\subsection{Additional definitions and notations.}
 Recall that for a non-negative r.v. $X$ with finite mean $\E[X] > 0$, one can define the so-called residual life distribution of $X$, $R(X)$, as follows.  Namely, for all $z \geq 0$,
\begin{equation}\label{rezz}
\pr\big( R(X) > z \big) = (\E[X])^{-1} \int_{z}^{\infty} \pr(X > y) dy.
\end{equation}
Recall that associated with a non-negative r.v. $X$, an equilibrium renewal process with renewal distribution $X$ is a counting process in which the first inter-event time is distributed as $R(X)$, and all subsequent inter-event times are drawn i.i.d. distributed as $X$; an ordinary renewal process with renewal distribution $X$ is a counting process in which all inter-event times are drawn i.i.d. distributed as $X$.  We will consider a FCFS $GI/GI/n$ queue with inter-arrival distribution $U$ and processing time distribution $V$, where we suppose that $U$ and $V$ are non-negative, have finite mean, and that $\pr(U = 0) = \pr(V = 0) = 0$.  Let $\lbrace {\mathcal N}_i, i \in [1,n] \rbrace$ denote a set of $n$ i.i.d. equilibrium renewal processes with renewal distribution $V$.  Let ${\mathcal A}$ denote an independent equilibrium renewal process with renewal distribution $U$.  For $t \geq 0$ we let $N_i(t)$ denote the number of renewals in ${\mathcal N}_i$ on $[0,t]$; and for $t_1 \leq t_2 \in \bbr^+$, we let $N_i(t_1,t_2)$ denote the number of renewals in ${\mathcal N}_i$ on $[t_1,t_2]$.  We define $A(t)$ and $A(t_1,t_2)$ analogously to describe the number of renewals in ${\mathcal A}$ over various increments.  We w.l.o.g. suppose that ${\mathcal A}$ and $\lbrace {\mathcal N}_i, i \in [1,n] \rbrace$ are constructed so as to be r.c.l.l. on $[0,\infty)$, where for $t \in \bbr^+$ we let $N_i(t^-)$, $A(t^-)$ denote the corresponding left limits.  For $t \geq 0$, we let $dA(t) \stackrel{\Delta}{=} A(t) - A(t^-)$ and $dN_i(t) \stackrel{\Delta}{=} N_i(t) - N_i(t^-)$.  Also, for $i \in [1,n]$, let $V^1_i$ denote the residual life of ${\mathcal N}_i$ at time 0.  Thus $V^1_i$ is distributed as $R(X)$ for all $i$.  For $j \geq 2$ and $i \in [1,n]$, let $V^j_i$ denote the length of the $(j-1)$st renewal interval to be initiated after time $0$ in ${\mathcal N}_i$.  Let $U^1$ and $U^j$ denote the analogous quantities for ${\mathcal A}$.  Also, for a set $E$, let $|E|$ denote its cardinality.
\\\indent Let ${\mathcal Q}$ denote the FCFS $GI/GI/n$ queue with initial conditions s.t. there is a single job with initial remaining processing time $V^1_i$ on server $i$, $i \in [1,n]$; the initial inter-arrival time is $U^1$; and there are zero jobs waiting in queue.  Also, let $Q(t)$ denote the total number in system at time $t$ in ${\mathcal Q}$.  
\subsection{Main upper bound result.}
We now establish an upper bound for $Q(t)$, by considering a modified queueing system which stochastically dominates the true system, and has a tractable steady-state distribution.  In \citet{GG.10c}, the authors considered a modified queueing system with arrival and service process nearly identical to that of ${\mathcal Q}$, which the exception that all servers were kept busy at all times by adding an artificial arrival whenever a server would otherwise have gone idle.  Although this was sufficient for analyzing the large deviations properties of the queue length, it will not suffice for our purposes, as in this modified system all servers are always busy, precluding the derivation of any meaningful bounds for the s.s.p.d. or number of idle servers.  To remedy this, we will consider a more sophisticated upper-bounding system, which considerably extends the ideas of \citet{GG.10c}.   Using this modified system, we will prove the following bounds for $Q(t)$.
\ \\\\\begin{theorem}\label{4ubound1}
For all $t,x \geq 0$, $\pr\big( Q(t) \geq x \big)$ is at most
\begin{eqnarray*}
&\ &\ \ \mathop{\inf_{\delta \in [0,t]}}_{\eta \in [0,n]}
\pr\Bigg( 
\max \bigg(\sup_{0 \leq s \leq \delta} \big( A(s) - \sum_{i=1}^{\eta} N_i(s) \big)\ \ ,\ \ 
\sup_{\delta \leq s \leq t} \big( A(s) - \sum_{i=1}^n N_i(s) \big) + \sum_{i = \eta + 1}^n N_i(\delta) \bigg)
\\&\ &\ \ \ \ \ \ \ \ \ \ \ \ \ \ + \sum_{i = \eta + 1}^n I(N_i(\delta) = 0)\ \ \geq\ \ x -\eta\Bigg).
\end{eqnarray*}
If in addition $Q(t)$ converges weakly to a steady-state distribution $Q(\infty)$ as $t \rightarrow \infty$, then for all $x \geq 0$, $\pr\big( Q(\infty) \geq x \big)$ is at most
\begin{eqnarray*}
&\ &\ \ \mathop{\inf_{\delta \geq 0}}_{\eta \in [0,n]}
\pr\Bigg( 
\max \bigg(\sup_{0 \leq t \leq \delta} \big( A(t) - \sum_{i=1}^{\eta} N_i(t) \big)\ \ ,\ \ \sup_{t \geq \delta} \big( A(t) - \sum_{i=1}^n N_i(t) \big) + \sum_{i = \eta + 1}^n N_i(\delta) \bigg)
\\&\ &\ \ \ \ \ \ \ \ \ \ \ \ \ \ + \sum_{i = \eta + 1}^n I(N_i(\delta) = 0)\ \ \geq\ \ x -\eta\Bigg).
\end{eqnarray*}
\end{theorem}
We note that the bounds of \citet{GG.10c} can be recovered by setting $\delta = t, \eta = n$ in the first part of the theorem.  We note that (in general) the choice of $\delta = t$ leads to simplified expressions (even for non-trivial choices of $\eta$), and we will later present an asymptotic (under the H-W scaling) theorem along these lines (see Corollary\ \ref{mainasymptotic2}).
\subsubsection{Outline of proof of Theorem\ \ref{4ubound1}.}
In the remainder of this section, we complete the proof of Theorem\ \ref{4ubound1}.  In Section\ \ref{newusystem}, we define our novel bounding system, first providing an informal description (to build intuition) in Section\ \ref{informal1}, and then a formal construction in Section\ \ref{formal1}.  In Section\ \ref{analyzeus} we analyze our bounding system, proving that its distribution obeys a certain Lindley-type recursion and can be expressed in terms of the suprema of certain one-dimensional random walks.  In Section\ \ref{compus}, we prove a relevant stochastic comparison result, and explicitly construct the original queueing system and our bounding system on a common probability space so that they satisfy this comparison result.  Finally, in Section\ \ref{proveit1}, we combine these results to complete the proof of Theorem\ \ref{4ubound1}.
\subsection{Novel bounding system: extra arrivals and servers that break down.}\label{newusystem}
In this section, we describe a novel bounding system $\overline{\mathcal Q}_{\eta,\gamma}$  for ${\mathcal Q}$, which we will use to prove our main results.
\subsubsection{Heuristic description.}\label{informal1}
Before providing a formal construction of our bounding system, we begin with a heuristic description to build intuition.  Suppose that we fix some integer $\eta \in [0,n]$ and real number $\gamma \in [0,t]$, and let $\overline{\mathcal Q}_{\eta,\gamma}$ be the following queueing system.  On $[0,\gamma]$, $\overline{\mathcal Q}_{\eta,\gamma}$ is identical to the bounding system considered in \citet{GG.10c}, namely that in which an artificial arrival is added whenever a server would otherwise have gone idle.  However, in $\overline{\mathcal Q}_{\eta,\gamma}$, at time $\gamma$ servers with index $i \in [\eta + 1,n]$ ``break down", in the sense that after completing the jobs which they are already processing at time $\gamma$, they cannot process any further jobs.  To stay within the framework of $G/G/n$ queues (and as it will be convenient in our analysis), we formally implement this ``server breakdown" by allowing those severs to continue processing jobs, but adding an artificial arrival any time any of those servers completes any jobs on $[\gamma,\infty)$ (irregardless of whether that service completion would have caused idling or not), effectively ``neutralizing" the ability of those servers to process jobs.  Furthermore, on $[\gamma,t]$, servers with index $i \in [1,\eta]$ are kept busy at all times by adding an artificial arrival whenever one of these $\eta$ servers would otherwise have gone idle.  We will prove that modulo the jobs being processed on servers with index $i \in [\eta + 1,n]$ at time $\gamma$, which result in a correction factor of the form $\sum_{i = \eta + 1}^n I(N_i(\delta) = 0)$, we can bound $Q(t)$ by the number of jobs on servers with index $i \in [1,\eta]$, plus the number of jobs waiting in queue, in this modified system.  Since in general $\eta \leq n$, this new modified system can have strictly less than $n$ jobs (if e.g. $\eta < n$ and there are 0 jobs waiting in queue), and we can get non-trivial bounds for the s.s.p.d. by computing the probability that in this modified system the total number of jobs on servers with index $i \in [1,\eta]$, plus the number of jobs waiting in queue, is less than $n$.  Also, by applying this technique with different values of $\eta$ and $\gamma$, we will be able to treat the different regimes of interest (e.g. $B \rightarrow \infty$ and $B \rightarrow 0$).
\subsubsection{Formal construction.}\label{formal1} We now provide a formal construction for $\overline{\mathcal Q}_{\eta,\gamma}$, building on the simpler construction provided in \citet{GG.10c}.  We proceed by defining processes $\lbrace \overline{A}_{\eta,\gamma}(t), t \geq 0 \rbrace$ and $\lbrace \overline{q}_{\eta,\gamma}(t), t \geq 0 \rbrace$, which will become the arrival and queue-length processes, respectively, for our bounding system.  For notational simplicity, in our construction we suppose $\gamma > 0$, as the case $\gamma = 0$ will then follow by continuity.  Let $\tau_{0} \stackrel{\Delta}{=} 0$, and let $\lbrace \tau_{k}, k \geq 1 \rbrace$ denote the ordered sequence of event times in the pooled renewal process $\lbrace A(t) + \sum_{i=1}^{n} N_i(t), t \geq 0 \rbrace$.  Note that we may w.l.o.g. condition on the event that $dA(\tau_{k}) + \sum_{i=1}^{n} dN_i(\tau_{k}) = 1$ for all $k \geq 1$, as this occurs w.p.1 since $R(S)$ and $R(A)$ are continuous r.v.s, and by assumption $\pr(A = 0) = \pr(S = 0) = 0$.  For clarity of exposition, we also implicitly condition on the occurence of other probability-one events, e.g. that for any fixed $\gamma$ no relevant events occur exactly at time $\gamma$.  Let $\overline{A}_{\eta,\gamma}(t) = 0$ and $\overline{q}_{\eta,\gamma}(t) = n$ for all $t \in [0, \tau_1)$.  Let $k^{\gamma}$ denote $\min\lbrace j \in Z^+ : \tau_j \geq \gamma \rbrace$, and set $\tau^{\gamma} \stackrel{\Delta}{=} \tau_{k^{\gamma}}$.  Now, we define the processes of interest inductively on $[\tau_1,\tau^{\gamma})$.  Suppose that for some $j \in [1, k^{\gamma} - 1]$, we have defined $\overline{A}_{\eta,\gamma}(t)$ and $\overline{q}_{\eta,\gamma}(t)$ for all $t \in [0, \tau_{j})$.  We now define these processes on $\big[\tau_{j}, \tau_{j+1})$.   Let
\begin{eqnarray}
\overline{A}_{\eta,\gamma}(\tau_{j}) &\stackrel{\Delta}{=} \begin{cases} 
\overline{A}_{\eta,\gamma}(\tau_{j}^-) + 1 & \text{if}\ dA(\tau_{j}) = 1;\\ 
\overline{A}_{\eta,\gamma}(\tau_{j}^-) + 1 & \text{if}\ \sum_{i=1}^{n} dN_i(\tau_{j}) = 1\ \text{and}\ \overline{q}_{\eta,\gamma}(\tau_{j}^-) \leq n;\\ 
\overline{A}_{\eta,\gamma}(\tau_{j}^-) & \text{otherwise}.
\end{cases}
\label{dddef1}
\end{eqnarray}
Similarly, let
\begin{eqnarray}
\overline{q}_{\eta,\gamma}(\tau_{j}) &\stackrel{\Delta}{=} \begin{cases} 
\overline{q}_{\eta,\gamma}(\tau_{j}^-) + 1 & \text{if}\ dA(\tau_{j}) = 1;\\ 
\overline{q}_{\eta,\gamma}(\tau_{j}^-) & \text{if}\ \sum_{i=1}^{n} dN_i(\tau_{j}) = 1\ \text{and}\ \overline{q}_{\eta,\gamma}(\tau_{j}^-) \leq n;\\ 
\overline{q}_{\eta,\gamma}(\tau_{j}^-) - 1& \text{otherwise}.
\end{cases}
\label{dddef2}
\end{eqnarray}
Also, set $\overline{A}_{\eta,\gamma}(t) = \overline{A}_{\eta,\gamma}(\tau_j)$ and $\overline{q}_{\eta,\gamma}(t) = \overline{q}_{\eta,\gamma}(\tau_j)$ for all $t \in (\tau_j, \tau_{j+1})$.  This completes our definition of these processes on $[0, \tau^{\gamma})$.  We now define these processes on $[\tau^{\gamma},\infty)$, again proceeding by induction.  Suppose that for some $j \geq k^{\gamma}$, we have defined $\overline{A}_{\eta,\gamma}(t)$ and $\overline{q}_{\eta,\gamma}(t)$ for all $t \in [0, \tau_{j})$.  We now define these processes on $\big[\tau_{j}, \tau_{j+1})$.   Let
\begin{eqnarray}
\overline{A}_{\eta,\gamma}(\tau_{j}) &\stackrel{\Delta}{=} \begin{cases} 
\overline{A}_{\eta,\gamma}(\tau_{j}^-) + 1 & \text{if}\ dA(\tau_{j}) = 1;\\ 
\overline{A}_{\eta,\gamma}(\tau_{j}^-) + 1 & \text{if}\ \sum_{i=1}^{\eta} dN_i(\tau_j) = 1\ \text{and}\ \overline{q}_{\eta,\gamma}(\tau_{j}^-) \leq n;\\ 
\overline{A}_{\eta,\gamma}(\tau_{j}^-) + 1 & \text{if}\ \sum_{i= \eta + 1}^{n} dN_i(\tau_j) = 1;\\
\overline{A}_{\eta,\gamma}(\tau_j^-) & \text{otherwise}.
\end{cases}
\label{dddef1}
\end{eqnarray}
Similarly, we define
\begin{eqnarray}
\overline{q}_{\eta,\gamma}(\tau_{j}) &\stackrel{\Delta}{=} \begin{cases} 
\overline{q}_{\eta,\gamma}(\tau_{j}^-) + 1 & \text{if}\ dA(\tau_{j}) = 1;\\ 
\overline{q}_{\eta,\gamma}(\tau_{j}^-) & \text{if}\ \sum_{i=1}^{\eta} dN_i(\tau_j) = 1\ \text{and}\ \overline{q}_{\eta,\gamma}(\tau_{j}^-) \leq n;\\ 
\overline{q}_{\eta,\gamma}(\tau_{j}^-) & \text{if}\ \sum_{i= \eta + 1}^{n} dN_i(\tau_{j}) = 1;\\ 
\overline{q}_{\eta,\gamma}(\tau_{j}^-) - 1& \text{otherwise}.
\end{cases}
\label{dddef3}
\end{eqnarray}
Also, set $\overline{A}_{\eta,\gamma}(t) = \overline{A}_{\eta,\gamma}(\tau_j)$ and $\overline{q}_{\eta,\gamma}(t) = \overline{q}_{\eta,\gamma}(\tau_j)$ for all $t \in (\tau_j, \tau_{j+1})$.  This completes our definition of these processes on $[\tau^{\gamma},\infty)$.
\\\indent Now, we use these processes to formally define $\overline{\mathcal Q}_{\eta,\gamma}$ as the following $G/G/n$ queue.  The arrival process is $\lbrace \overline{A}_{\eta,\gamma}(t), t \geq 0 \rbrace$.  We assign processing times as follows.  For $i \in [1,n]$, server $i$ initially has a single job with initial remaining processing time $V^1_i$; and the $j$th job assigned to server $i$ (after time 0) receives processing time $V^{j+1}_{i}$.  Note that processing times are thus i.i.d. distributed as $V$, with all initial remaining processing times distributed as $R(V)$.  As it will be convenient in analyzing the couplings used in our analysis, we will have certain arriving jobs join the front (as opposed to the back) of the queue, as follows.  On $[0,\gamma]$, all arriving jobs join the back of the queue.  On $[\gamma,\infty)$, any job which arrives at a time belonging to $\lbrace t \geq \gamma : \sum_{i = \eta + 1}^n dN_i(t) = 1 \rbrace$ joins the front of the queue; all other jobs which arrive on $[\gamma,\infty)$ join the back of the queue.
This deviation from FCFS will allow us to easily account for the ``extra arrivals" triggered by departures on servers with index $i \in [\eta + 1,n]$ on $[\gamma,\infty)$.  Furthermore, any time a job completes service, if the queue is strictly positive the next job to begin service is always the job at the front of the queue (i.e. standard head-of-line service discipline).  In the case that there are multiple empty servers, suppose a job always begins service on the server with lowest index.  This completely specifies the description of $\overline{\mathcal Q}_{\eta,\gamma}$ as a $G/G/n$ queue.
\subsection{Analysis of $\overline{\mathcal Q}_{\eta,\gamma}$.}\label{analyzeus}
In this section we analyze $\overline{\mathcal Q}_{\eta,\gamma}$ and prove that the system has several properties which we will need for our stochastic comparison arguments.  We begin with an important lemma and associated corollary, which characterize the dynamics of 
$\overline{\mathcal Q}_{\eta,\gamma}$, proving e.g. that no server ever idles, and that $\overline{Q}_{\eta,\gamma}(t) = \overline{q}_{\eta,\gamma}(t).$  These results follow from arguments nearly identical to those used to prove very similar results for a closely related bounding system in \citet{GG.10c}, i.e. in the proofs of Lemma 1 and Corollary 1 of that paper.  As such, we omit the proofs, instead referring the interested reader to \citet{GG.10c} for details.  For $i \in [1,n]$, let $V_i(t)$ denote the residual life (i.e. time until the next renewal) in ${\mathcal N}_i$ at time $t$, and $U(t)$ denote the analogous quantity for ${\mathcal A}$.  

\begin{lemma}\label{4busybeaver2}
For the $G/G/n$ queue $\overline{\mathcal Q}_{\eta,\gamma}$, the following are true.
\begin{enumerate}[(i)]
\item For $i \in [1,n]$, exactly one job departs from server $i$ at each time $t \in \lbrace \sum_{l=1}^j V^l_i, j \geq 1 \rbrace$, and there are no other departures from server $i$.  
\item No server ever idles, and the remaining processing time of the job on server $i$ (at time $t$) equals $V_i(t)$.
\item For all $t \geq 0$, $\overline{Q}_{\eta,\gamma}(t) = \overline{q}_{\eta,\gamma}(t).$  
\item \label{bbrecurse1}
For all $j \in [0, k^{\gamma} - 1]$, 
$$\overline{Q}_{\eta,\gamma}(\tau_{j + 1}) - n = \max\big(0, \overline{Q}_{\eta,\gamma}(\tau_{j}) - n + d A(\tau_{j+1}) - \sum_{i=1}^n d N_i(\tau_{j+1}) \big).$$
\item \label{bbrecurse2} For all $j \geq k^{\gamma}$, 
$$\overline{Q}_{\eta,\gamma}(\tau_{j + 1}) - n = \max\big(0, \overline{Q}_{\eta,\gamma}(\tau_{j}) - n + d A(\tau_{j+1}) - \sum_{i=1}^{\eta} d N_i(\tau_{j+1}) \big).$$
\end{enumerate}
\end{lemma}

As in \citet{GG.10c}, we may then ``unfold" the recursions indicated in Lemma\ \ref{4busybeaver2}.(\ref{bbrecurse1}) - (\ref{bbrecurse2}) to derive the following representation for $\overline{Q}_{\eta,\gamma}(t)$.  For a further discussion of such max-plus recursions and their connection to e.g. Skorohod problems, we refer the reader to the discussion in \citet{GG.10c}.

\begin{corollary}\label{4busybeaver2cor}
$$\overline{Q}_{\eta,\gamma}(\gamma) - n = \sup_{0 \leq u \leq \gamma} \big( A(\gamma - u, \gamma) - \sum_{i=1}^{n} N_i(\gamma - u, \gamma) \big),$$ and for all $t \geq \gamma$,
$\overline{Q}_{\eta,\gamma}(t) - n$ equals
$$
\max\bigg( \sup_{0 \leq u \leq t - \gamma} \big( A(t - u,t) - \sum_{i=1}^{\eta} N_i(t - u,t) \big)\ \ ,\ \  \sup_{0 \leq u \leq \gamma} \big( A(\gamma - u, \gamma) - \sum_{i=1}^{n} N_i(\gamma - u, \gamma) \big) + A(\gamma, t) - \sum_{i=1}^{\eta} N_i(\gamma, t) \bigg).$$
\end{corollary}

\subsection{Stochastic comparison result and using $\overline{\mathcal Q}_{\eta,\gamma}$ to bound ${\mathcal Q}$.}\label{compus}

In this section we state a novel stochastic comparison result for queues, and use this result to prove that $\overline{Q}_{\eta,\gamma}(t)$ yields an upper bound for $Q(t)$.  The result follows from arguments nearly identical to those used to prove a very similar result in \citet{GG.10c}, i.e. in the proof of Lemma 2 of that paper.  As such, we omit the proofs, instead referring the interested reader to \citet{GG.10c} for details.

\begin{lemma}\label{qom}
Let ${\mathcal Q}^1$ and ${\mathcal Q}^2$ be two $G/G/n$ queues with finite, strictly positive inter-arrival and processing times, and $Q^1(t), Q^2(t)$ the corresponding number in system at time $t$.  Let $\lbrace T^i_k, k \geq 1 \rbrace$ denote the ordered sequence of arrival times to ${\mathcal Q}^i$, $i \in \lbrace 1,2 \rbrace$.  Let $S^i_k$ denote the processing time assigned to the job that arrives to ${\mathcal Q}^i$ at time $T^i_k$, $k \geq 1, i \in \lbrace 1,2 \rbrace$.  Also, although we allow arriving jobs to join different locations in the queue (i.e. some arriving jobs may join the front of the queue instead of the back of the queue), we assume both queues operate under the standard head-of-line service discipline.  Further suppose that:
\begin{enumerate}[(i)]
\item The initial number in system in both ${\mathcal Q}^1$ and ${\mathcal Q}^2$ is $n$, and the initial remaining processing time of the job on server $i$ is the same in both ${\mathcal Q}^1$ and ${\mathcal Q}^2$ for all $i \in [1,n]$;\label{as1}
\item $\lbrace T^1_k, k \geq 1\rbrace$ is a subsequence of $\lbrace T^2_k, k \geq 1 \rbrace$;\label{as3}
\item For all $k \geq 1$, the job that arrives to ${\mathcal Q}^1$ (${\mathcal Q}^2$) at time $T^1_k$ joins the back of the queue in ${\mathcal Q}^1$ (${\mathcal Q}^2$);\label{as3a}
\item For all $k \geq 1$, the job that arrives to ${\mathcal Q}^2$ at time $T^1_k$ is assigned processing time $S^1_k$, the same processing time assigned to the job which arrives to ${\mathcal Q}^1$ at that time.\label{as4}
\end{enumerate}
\ \\For $i \in \lbrace 1,2 \rbrace$, let $D^i_k$ denote the time at which the job which arrives to ${\mathcal Q}^i$ at times $T^1_k$ departs from ${\mathcal Q}^i$.
Let $D(t)$ denote the set of jobs which arrive to ${\mathcal Q}^2$ at a time belonging to $[0,t) \setminus \lbrace T^1_k, k \geq 1 \rbrace$, and which have not yet departed from ${\mathcal Q}^2$ by time $t$.  Namely, $D(t)$ is the set of jobs which are in ${\mathcal Q}^2$ at time $t$ and which had no corresponding arrival in ${\mathcal Q}^1$.
\\\\Then for all $k \geq 1$, $D^2_k \geq D^1_k$.  Also, $Q^2(t) \geq Q^1(t) + |D(t)|$ for all $t \geq 0$.
\end{lemma}

We now apply Lemma\ \ref{qom} to $\overline{\mathcal Q}_{\eta,\gamma}$ and ${\mathcal Q}$. 
\begin{corollary}\label{qomcor}
For all $\eta \in [0,n]$ and $\gamma > 0$, one may construct $\overline{\mathcal Q}_{\eta,\gamma}$ and ${\mathcal Q}$ on a common probability space s.t. w.p.1, for all $t \geq \gamma$, 
$\overline{Q}_{\eta,\gamma}(t) - n + \eta + \sum_{i = \eta + 1}^n I\big( V_i(\gamma) > t - \gamma \big) \geq Q(t)$.
\end{corollary}
\proof{Proof:}
We construct $\overline{\mathcal Q}_{\eta,\gamma}$ and ${\mathcal Q}$ on the same probability space.  We assign ${\mathcal Q}$ and $\tilde{\mathcal Q}$ the same initial conditions, and let $\lbrace A(t), t \geq 0 \rbrace$ be the arrival process to ${\mathcal Q}$ on $(0,\infty)$.  Let $\lbrace t_k , k \geq 1 \rbrace$ denote the ordered sequence of event times in $\lbrace A(t), t \geq 0 \rbrace$.  It follows from the construction of $\lbrace \overline{A}_{\eta,\gamma}(t), t \geq 0 \rbrace$ that $\lbrace t_k , k \geq 1 \rbrace$ is a subsequence of the set of event times in $\lbrace \overline{A}_{\eta,\gamma}(t), t \geq 0 \rbrace$.  We let the processing time assigned to the arrival to ${\mathcal Q}$ at time $t_k$ equal the processing time assigned to the arrival to $\overline{\mathcal Q}_{\eta,\gamma}$ at time $t_k$, $k \geq 1$.  It follows that w.p.1 ${\mathcal Q}$ and $\overline{\mathcal Q}_{\eta,\gamma}$ satisfy the conditions of Lemma\ \ref{qom}.  Furthermore, we may conclude from Lemma\ \ref{4busybeaver2} and a straightforward induction that (w.p.1) for $i \in [\eta + 1,n]$, at every time $t \in \lbrace \sum_{l=1}^j V^l_i, j \geq 1 \rbrace \bigcap [\gamma,\infty)$, there is an arrival to $\overline{\mathcal Q}_{\eta,\gamma}$, and no arrival to ${\mathcal Q}$.  Also, the corresponding arrival to $\overline{\mathcal Q}_{\eta,\gamma}$ immediately begins processing on server $i$ (by virtue of the corresponding job arriving to the front of the queue).  It thus follows from a straightforward induction that for all $i \in [\eta + 1, n]$ and times $t \geq \gamma + V_i(\gamma)$, in $\overline{\mathcal Q}_{\eta,\gamma}$ server $i$ is occupied by a job which had no corresponding arrival in ${\mathcal Q}$.  This implies that (in the terminology of Lemma\ \ref{qom}) $|D(t)| \geq \sum_{i = \eta + 1}^n I\big( V_i(\gamma) \leq t - \gamma \big)$ for all $t \geq \gamma$.  Combining with Lemma\ \ref{qom}, and the fact that 
$- \sum_{i = \eta + 1}^n I\big( V_i(\gamma) \leq t - \gamma \big) = - (n - \eta) + \sum_{i = \eta + 1}^n I\big( V_i(\gamma) > t - \gamma \big)$, completes the proof.  $\Halmos$.
\endproof

\subsection{Proof of Theorem\ \ref{4ubound1}.}\label{proveit1}
With Corollaries\ \ref{4busybeaver2cor}\ and\ \ref{qomcor} in hand, we now complete the proof of Theorem\ \ref{4ubound1}.
\proof{Proof:}[Proof of Theorem\ \ref{4ubound1}]
Combining Corollaries\ \ref{4busybeaver2cor}\ and\ \ref{qomcor}, we conclude that for all $\eta \in [0,n], \gamma > 0$, $t \geq \gamma$, and $x \geq 0$, $\pr\big( Q(t) \geq x \big)$ is at most
\begin{eqnarray}
&\ &\ \pr\Bigg(\max\bigg( \sup_{0 \leq u \leq t - \gamma} \big( A(t - u,t) - \sum_{i=1}^{\eta} N_i(t - u,t) \big)\ \ , \label{long1a}
\\&\ &\ \ \ \ \ \ \ \ \ \ \ \ \ \  \sup_{0 \leq u \leq \gamma} \big( A(\gamma - u, \gamma) - \sum_{i=1}^{n} N_i(\gamma - u,\gamma) \big) + A(\gamma, t) - \sum_{i=1}^{\eta} N_i(\gamma, t) \bigg) \nonumber
\\&\ &\ \ \ \ \ \ \ \ \ \ \ \ \ \  + \sum_{i = \eta + 1}^n I\big( V_i(\gamma) > t - \gamma \big) \geq x - \eta \Bigg). \nonumber
\end{eqnarray}
\endproof
By elementary renewal theory, and specifically the well-known properties of equilibrium renewal processes (see \citet{Cox.70}), for all $t \geq \gamma \geq 0$, the joint distribution of 
$$A\big(t - u, t\big)_{0 \leq u \leq t}\ \ ,\ \ A\big(\gamma - u, \gamma\big)_{0 \leq u \leq \gamma}\ \ ,\ \ \sum_{i=1}^n N_i\big(t - u, t\big)_{0 \leq u \leq t}\ \ ,\ \ \sum_{i=1}^n N_i\big(\gamma - u, \gamma\big)_{0 \leq u \leq \gamma}, 
\sum_{i = \eta + 1}^n I\big( V_i(\gamma) > t - \gamma \big)$$
is the same as that of 
$$A(u)_{0 \leq u \leq t}\ \ ,\ \ A\big(t - \gamma, t - \gamma + u\big)_{0 \leq u \leq \gamma}\ \ ,\ \ \sum_{i=1}^n N_i(u)_{0 \leq u \leq t}\ \ ,\ \ \sum_{i=1}^n N_i\big(t - \gamma, t - \gamma + u\big)_{0 \leq u \leq \gamma}\ \ ,\ \ \sum_{i = \eta + 1}^n I\big( V^1_i > t - \gamma \big).$$
Combining with (\ref{long1a}), we conclude that for all $\eta \in [0,n], \gamma > 0$, $t \geq \gamma$, and $x \geq 0$, $\pr\big( Q(t) \geq x \big)$ is at most
\begin{eqnarray}
&\ &\ \pr\Bigg(\max\bigg( \sup_{0 \leq u \leq t - \gamma} \big( A(u) - \sum_{i=1}^{\eta} N_i(u) \big)\ \ , \label{ddd1}
\\&\ &\ \ \ \ \ \ \ \ \ \ \ \ \ \  \sup_{0 \leq u \leq \gamma} \big( A(t - \gamma , t - \gamma + u) - \sum_{i=1}^{n} N_i(t - \gamma, t - \gamma + u) \big) + A(t - \gamma) - \sum_{i=1}^{\eta} N_i(t - \gamma) \bigg) \nonumber
\\&\ &\ \ \ \ \ \ \ \ \ \ \ \ \ \  + \sum_{i = \eta + 1}^n I\big( V^1_i > t - \gamma \big) \geq x - \eta \Bigg). \nonumber
\end{eqnarray}
Letting $\delta \stackrel{\Delta}{=} t - \gamma$, note that
\begin{eqnarray}
\ &\ &\ \sup_{0 \leq u \leq \gamma} \big( A(t - \gamma , t - \gamma + u) - \sum_{i=1}^{n} N_i(t - \gamma, t - \gamma + u) \big) + A(t - \gamma) - \sum_{i=1}^{\eta} N_i(t - \gamma) \nonumber
\\&\ &\ \ \ =\ \ \ \sup_{0 \leq u \leq t - \delta} \big( A(\delta , \delta + u) - \sum_{i=1}^{n} N_i(\delta, \delta + u) \big) + A(\delta) - \sum_{i=1}^{\eta} N_i(\delta) \nonumber
\\&\ &\ \ \ =\ \ \ \sup_{0 \leq u \leq t - \delta} \big( A(\delta + u) - \sum_{i=1}^{n} N_i(\delta + u) \big) + \sum_{i = \eta + 1}^n N_i(\delta) \nonumber
\\&\ &\ \ \ =\ \ \ \sup_{\delta \leq u \leq t} \big( A(u) - \sum_{i=1}^{n} N_i(u) \big) + \sum_{i = \eta + 1}^n N_i(\delta). \label{ddd2}
\end{eqnarray}
Combining (\ref{ddd1}) and (\ref{ddd2}) with the fact that $I\big( V^1_i > \delta \big) = I\big(N_i(\delta) = 0\big)$, we conclude that for all $t > 0, \eta \in [0,n], \delta \in (0,t)$, and $x \geq 0$, $\pr\big( Q(t) \geq x \big)$ is at most
\begin{eqnarray}
&\ &\ \pr\Bigg(\max\bigg( \sup_{0 \leq u \leq \delta} \big( A(u) - \sum_{i=1}^{\eta} N_i(u) \big)\ \ ,\ \ \sup_{\delta \leq u \leq t} \big( A(u) - \sum_{i=1}^{n} N_i(u) \big) + \sum_{i = \eta + 1}^n N_i(\delta) \bigg) \nonumber
\\&\ &\ \ \ \ \ \ \ \ \ \ \ \ \ \  + \sum_{i = \eta + 1}^n I\big(N_i(\delta) = 0\big) \geq x - \eta \Bigg). \nonumber
\end{eqnarray}
Minimizing over $\eta$ and $\delta$ completes the proof of the first part of the theorem.  Combining with the monotonicity of the supremum operator, the continuity of probability measures, and the definition of weak convergence completes the proof of the second part of the theorem and overall result.  $\Halmos$.
\endproof

\section{Asymptotic bound in the H-W regime.}\label{4upperasymptoticsec}
In this section, we use Theorem\ \ref{4ubound1} to bound the FCFS $GI/GI/n$ queue in the H-W regime, by proving an asymptotic analogue of Theorem\ \ref{4ubound1}.  
We note that in the remainder of the paper, we restrict our analysis to the steady-state queue length, as this is the primary motivation of our investigations.  However, we note that corresponding finite time horizon asymptotic bounds could be analogously derived using nearly identical arguments.  Suppose that the H-W and $T_0$ assumptions hold.  
For a r.v. $X$, let $V[X]$ denote the variance of $X$.  For r.v.s $X,Y$, let $V[X,Y] \stackrel{\Delta}{=} \E[XY]-\E[X]\E[Y]$ denote the covariance of $X$ and $Y$.  Recall that a Gaussian process on $\reals$ is a stochastic process $Z(t)_{t \geq 0}$ s.t. for any finite set of times $t_1,\ldots,t_k$, the vector $\big( Z(t_1),\ldots,Z(t_k) \big)$ has a Gaussian distribution.  A Gaussian process $Z$ is known to be uniquely determined by its mean function $\E[Z(t)]$ and covariance function $V[Z(s) , Z(t)]$, and we refer the reader to \citet{D.44}, \citet{IR.78}, \citet{A.90}, \citet{MR.06},  and the references therein for details on existence, continuity, etc.  It is proven by \citet{Whitt.85} Theorem 2 that there exists a continuous Gaussian process ${\mathcal D}$ s.t. $\E[{\mathcal D}(t)] = 0, V[{\mathcal D}(s) , {\mathcal D}(t)] = \E[\big( N_1(s) - \mu  s \big) \big( N_1(t) - \mu  t \big)]$ for all $s,t \geq 0$.  Let ${\mathcal A}$ denote the continuous Gaussian process s.t. $\E[{\mathcal A}(t)] = 0, V[{\mathcal A}(s) , {\mathcal A}(t)] = \mu c^2_A \min(s,t)$, namely ${\mathcal A}$ is a driftless Brownian motion.  Let ${\mathcal Z} \stackrel{\Delta}{=} {\mathcal A} - {\mathcal D}$, where ${\mathcal A}$ and ${\mathcal D}$ are independent.  For $\delta,\eta \geq 0$, $B > 0$, and $x \in \bbr$, let us define the event
$${\mathcal E}^{\delta,\eta}_{B,\infty}(x) \stackrel{\Delta}{=} 
\Bigg\lbrace \max \bigg( \sup_{0 \leq t \leq \delta} \big( {\mathcal Z}(t) + (\eta - B) \mu t \big)\ \ ,\ \ \sup_{t \geq \delta} \big( {\mathcal Z}(t) - B \mu t \big)\ +\ \eta \mu \delta \bigg)
\geq x + \eta \pr\big( R(S) \leq \delta \big) \Bigg\rbrace.$$
Then our main asymptotic upper bound is the following.
\begin{theorem}\label{mainasymptotic}
For all $B > 0$ and $x \in \reals$ , 
$$\limsup_{n \rightarrow \infty} \pr\bigg( \big( Q^n_B(\infty) - n \big)n^{-\frac{1}{2}} \geq x \bigg) \leq 
\mathop{\inf_{\delta \geq 0}}_{\eta \geq 0}
\pr\big( {\mathcal E}^{\delta,\eta}_{B,\infty}(x) \big).$$
\end{theorem}
We prove Theorem\ \ref{mainasymptotic} by first rewriting Theorem\ \ref{4ubound1} in a manner more amenable to asymptotic analysis, and then carrying out such an asymptotic analysis using the framework of weak convergence.  Let ${\mathcal A}^n_B$ denote an independent equilibrium renewal process with renewal distribution $A \lambda_{n,B}^{-1}$, $\lbrace {\mathcal N}^S_i, i \in [1,n] \rbrace$ denote a set of $n$ i.i.d. equilibrium renewal processes with renewal distribution $S$, and $A^n_B(t), N^S_i(t)$ denote the corresponding notation for the number of renewals up to time $t$.  Similarly, let $Z^n_B(t) \stackrel{\Delta}{=} A^n_B(t) - \sum_{i=1}^n N^S_i(t)$.  Also, we define $f_n(x) \stackrel{\Delta}{=} \lfloor n - x n^{\frac{1}{2}} \rfloor$, and
\begin{eqnarray*}
{\mathcal E}^{\delta,\eta}_{B,n}(x) &\stackrel{\Delta}{=}&\Bigg\lbrace
\max \bigg( \sup_{0 \leq t \leq \delta} \big( Z^n_B(t) + \sum_{i=f_n(\eta)+1}^{n} N^S_i(t) \big)\ \ \ ,\ \ \ \sup_{t \geq \delta} Z^n_B(t) + \sum_{ i = f_n(\eta) + 1 }^n N^S_i(\delta) \bigg) 
\\&\indent&\ \ \ \ \ + \sum_{i = f_n(\eta) + 1}^n I\big(N^S_i(\delta) = 0\big)\ \ \geq\ \ f_n(-x) - f_n(\eta)\Bigg\rbrace.
\end{eqnarray*}
Then our rewriting of Theorem\ \ref{4ubound1} is as follows.
\begin{lemma}
For all $B > 0$, $x \in \reals$, and $n$ s.t. $\lambda_{n,B} > 0$,
$$\pr\bigg( \big( Q^n_B(\infty) - n \big)n^{-\frac{1}{2}} \geq x \bigg) \leq \mathop{\inf_{\delta \geq 0}}_{\eta \in [0,n^{\frac{1}{2}}]} \pr\big( {\mathcal E}^{\delta,\eta}_{B,n}(x) \big),
$$
where $\eta$ is no longer restricted to the integers and may be any real number in that interval.
\end{lemma}
\proof{Proof:} $f_n$, when restricted to domain $[0, n^{\frac{1}{2}}]$, is a surjective function with range $[0,n]$.  It follows that Theorem\ \ref{4ubound1} is equivalent to asserting that $\pr\bigg( \big( Q^n_B(\infty) - n \big)n^{-\frac{1}{2}} \geq x \bigg)$ is at most
\begin{eqnarray}
&\ &\ \ \mathop{\inf_{\delta \geq 0}}_{\eta \in [0, n^{\frac{1}{2}}]} \label{asym1}
\pr\Bigg( 
\max \bigg(\sup_{0 \leq t \leq \delta} \big( A^n_B(t) - \sum_{i=1}^{f_n(\eta)} N^S_i(t) \big)\ \ \ ,\ \ \ \sup_{t \geq \delta} Z^n_B(t) + \sum_{i = f_n(\eta) + 1}^n N^S_i(\delta) \bigg) \nonumber
\\&\ &\ \ \ \ \ \ \ \ \ \ \ \ \ \ + \sum_{i = f_n(\eta) + 1}^n I\big(N^S_i(\delta) = 0\big)\ \ \geq\ \ f_n(-x) - f_n(\eta)\Bigg). \nonumber
\end{eqnarray}
Noting that $A^n_B(t) - \sum_{i=1}^{f_n(\eta)} N^S_i(t) = Z^n_B(t) + \sum_{i=f_n(\eta)+1}^{n} N^S_i(t)$ completes the proof.  $\Halmos$.
\endproof

\subsection{Preliminary weak convergence results.}
Before completing the proof of Theorem\ \ref{mainasymptotic}, we establish some preliminary weak convergence results to aid in the analysis of $\pr\big({\mathcal E}^{\delta,\eta}_{B,n}(x)\big)$.  For an excellent review of weak convergence, and the associated spaces (e.g. $D[0,T]$) and topologies/metrics (e.g. uniform, $J_1$), we refer the reader to \citet{W.02}.  The following weak convergence results follows immediately from the results of \citet{GG.10c}, specifically Theorem 7, Lemmas 6 - 7, Corollary 3, and Equation 17 of that paper.
\begin{theorem}[\citet{GG.10c}]\label{oldandweak}
\ \begin{enumerate}[(i)]
\item \label{4renewalfclt2} For any fixed $T \in [0,\infty)$, the sequence of processes $\big\lbrace n^{-\frac{1}{2}} \big( \sum_{i=1}^n N^S_i(t) - n \mu  t \big)_{0 \leq t \leq T}, n \geq 1 \big\rbrace$ converges weakly to ${\mathcal D}(t)_{0 \leq t \leq T}$ in the space $D[0,T]$ under the $J_1$ topology.
\item \label{4renewalfclt3}
For any fixed $B > 0$ and $T \in [0,\infty)$, the sequence of processes $\big\lbrace n^{-\frac{1}{2}} Z^n_B(t)_{0 \leq t \leq T}, n \geq 1 \big\rbrace$ converges weakly to 
$\big( {\mathcal Z}(t) - B \mu  t \big)_{0 \leq t \leq T}$ in the space $D[0,T]$ under the $J_1$ topology.
\item \label{showsup1} For any fixed $B,x > 0$, 
$\limsup_{T \rightarrow \infty} \limsup_{n \rightarrow \infty} \pr\big( n^{-\frac{1}{2}} \sup_{t \geq T} Z^n_B(t) \geq x\big) = 0$.  It follows that the sequence of r.v.s $\lbrace n^{-\frac{1}{2}} \sup_{t \geq 0} Z^n_B(t) , n \geq 1 \rbrace$ is tight.
\item \label{varlimitplus}
$\lim_{t \rightarrow \infty} \E[\big( t^{-\frac{1}{2}}{\mathcal Z}(t) \big)^2] = \mu(c^2_A + c^2_S)$.  Also, for all $c > 0$, the r.v.s $\sup_{t \geq 0} \big( {\mathcal Z}(t) - c t \big)$ and $\sup_{t \geq 0} \big( {\mathcal D}(t) - c t \big)$ are a.s. finite.
\end{enumerate}
\end{theorem}
We now state several additional weak convergence results, also for use in our analysis.  These results follow directly from Theorem\ \ref{oldandweak}, the strong law of large numbers, and the basic properties and definitions of weak convergence, combined with some straightforward algebraic manipulations.  As such, we omit the associated proofs for clarity of exposition.
\begin{corollary}\label{oldweakcor}
\ \begin{enumerate}[(i)]
\item \label{showsup1bb}
For any fixed $B > 0$, $\delta \geq 0$, and $\epsilon > 0$, there exists $T \in (\delta,\infty)$ s.t.
$\limsup_{n \rightarrow \infty} \pr \bigg( \sup_{ t \geq \delta } Z^n_B(t) \neq \sup_{t \in [\delta,T]} Z^n_B(t) \bigg) < \epsilon.$
\item \label{fnetalimit}
For any fixed $\eta, x \in \bbr$, 
$\lim_{n \rightarrow \infty} n^{-\frac{1}{2}}\big(n - f_n(\eta)\big) = \eta$, and $\lim_{n \rightarrow \infty} n^{-\frac{1}{2}}\big(f_n(-x) - f_n(\eta)\big) = x + \eta.$
\item \label{errors23gonea}
For any fixed $\eta, \delta \geq 0$, the sequence of processes\ $\big\lbrace n^{-\frac{1}{2}} \sum_{ i = f_n(\eta) + 1 }^n N^S_i(t)_{0 \leq t \leq \delta} , n \geq 1 \big\rbrace$\ converges weakly to $\eta \mu t_{0 \leq t \leq \delta}$ in the space $D[0,T]$ under the $J_1$ topology.
\item \label{errors23goneb} 
For any fixed $\eta, \delta \geq 0$, the sequence of r.v.s\ $\big\lbrace n^{-\frac{1}{2}}  \sum_{i = f_n(\eta) + 1}^n I\big(N^S_i(\delta) = 0\big) , n \geq 1 \big\rbrace$\ converges weakly to the constant\ $\eta \pr\big( R(S) > \delta \big)$.
\end{enumerate}
\end{corollary}
\subsection{Proof of Theorem\ \ref{mainasymptotic}.}
We now complete the proof of Theorem\ \ref{mainasymptotic}.
\\\proof{Proof:}[Proof of Theorem\ \ref{mainasymptotic}]
Note that
\begin{eqnarray*}
{\mathcal E}^{\delta,\eta}_{B,n}(x) &=&\Bigg\lbrace
\max \bigg( \sup_{0 \leq t \leq \delta} \big( n^{-\frac{1}{2}} Z^n_B(t) + n^{-\frac{1}{2}} \sum_{i=f_n(\eta)+1}^{n} N^S_i(t) \big)\ \ \ ,\ \ \ \sup_{t \geq \delta} n^{-\frac{1}{2}} Z^n_B(t) + n^{-\frac{1}{2}} \sum_{ i = f_n(\eta) + 1 }^n N^S_i(\delta) \bigg) 
\\&\indent&\ \ \ \ \ + n^{-\frac{1}{2}} \sum_{i = f_n(\eta) + 1}^n I\big(N^S_i(\delta) = 0\big)\ \ \geq\ n^{-\frac{1}{2}} \big( f_n(-x) - f_n(\eta) \big) \Bigg\rbrace.
\end{eqnarray*}
Combining with Corollary\ \ref{oldweakcor}.(\ref{showsup1bb}) and a union bound, it follows that to prove the theorem, it suffices to demonstrate that for each fixed $B > 0$ and $\delta,\eta \geq 0$,  
\begin{eqnarray*}
&\ &\limsup_{T \rightarrow \infty}\ \ \limsup_{n \rightarrow \infty}\ \ \pr
\Bigg(
\max \bigg( \sup_{0 \leq t \leq \delta} \big( n^{-\frac{1}{2}} Z^n_B(t) + n^{-\frac{1}{2}} \sum_{i=f_n(\eta)+1}^{n} N^S_i(t) \big)\ \ \ ,
\\&\ &\ \ \ \ \ \ \ \ \ \ \ \ \ \ \ \ \ \ \ \ \ \ \ \ \ \ \ \ \ \ \ \ \ \ \ \ \sup_{t \in [\delta,T]} n^{-\frac{1}{2}} Z^n_B(t) + n^{-\frac{1}{2}} \sum_{ i = f_n(\eta) + 1 }^n N^S_i(\delta) \bigg)
\\&\ &\ \ \ \ \ \ \ \ \ \ \ \ \ \ \ \ \ \ \ \ \ \ \ \ \ \ \ \ +\ \ n^{-\frac{1}{2}} \sum_{i = f_n(\eta) + 1}^n I\big(N^S_i(\delta) = 0\big)\ \ \geq\ n^{-\frac{1}{2}} \big( f_n(-x) - f_n(\eta) \big) \Bigg) \ \ \ \leq\ \ \ \pr\big( {\mathcal E}^{\delta,\eta}_{B,\infty}(x)\big).
\end{eqnarray*}
The desired result then follows from Theorem\ \ref{oldandweak}, Corollary\ \ref{oldweakcor}, the continuity of the supremum map (see \citet{W.02}), and the Portmanteau Theorem (see \citet{B.99}).  $\Halmos$
\endproof
\section{Proof of bound for s.s.p.d. as $B \rightarrow \infty$.}\label{Binfinity}
In this section we complete the proof of Theorem\ \ref{largebtheorem}.  We begin by proving a modified variant of Theorem\ \ref{mainasymptotic}, which has the interpretation of setting $\delta = \infty$ in Theorem\ \ref{mainasymptotic}, and defer the proof to the appendix.
\begin{corollary}\label{mainasymptotic2}
For all $B > 0$ and $x \in \reals$ , $\limsup_{n \rightarrow \infty} \pr\bigg( \big( Q^n_B(\infty) - n \big)n^{-\frac{1}{2}} \geq x \bigg)$ is at most
$$
\inf_{\eta \in [0,B)} \pr\bigg( \sup_{t \geq 0} \big( {\mathcal Z}(t) - (B - \eta) \mu t \big)
\ \ \geq\ \ x + \eta \bigg).
$$
\end{corollary}
We next state a well-known result from the theory of Gaussian processes (see \citet{A.90} Inequality 2.4) which will be critical to our proof.
\begin{lemma}[\citet{A.90}]\label{intervalbound}
Let ${\mathcal X}$ denote any centered, continuous Gaussian process, and $T$ any bounded interval of $\reals^+$.  Let $\sigma^2_T \stackrel{\Delta}{=} \sup_{t \in T} \E[{\mathcal X}^2(t)]$, and suppose $\sigma^2_T < \infty$.  Then for all $\epsilon > 0$, there exists $M_{\epsilon}$, depending only on ${\mathcal X}, T$, and $\epsilon$, s.t. for all $x > M_{\epsilon}$,
$$\pr\big( \sup_{t \in T} {\mathcal X}(t) > x \big) \leq \exp\Bigg( - \bigg( (2 \sigma^2_T)^{-1} - \epsilon \bigg) x^2 \Bigg).$$
\end{lemma}
We will also need the following bounds on the variance of ${\mathcal Z}$, whose proofs we defer to the appendix (see Subsection\ \ref{proofzztop}).
\begin{lemma}\label{boundvarzztop}
For all $t \geq 0$, $V[{\mathcal D}(t)] \leq  \mu \big( 2 c^2_S + 1 \big) t$, and 
$V[{\mathcal Z}(t)] \leq \mu \big( 2 c^2_S + c^2_A + 1 \big) t.$
\end{lemma}
Finally, we will need the following bound on the Gaussian c.d.f., which follows from the results of \citet{abramowitz1999ia}.  Let $\Phi$ denote the cumulative distribution function of a standard normal distribution, $\phi$ the corresponding probability density function, and let us define $\Phi^c \stackrel{\Delta}{=} 1 - \Phi$.  Also, for $a \in \bbr$ and $b \geq 0$, let $N(a,b)$ denote a normally distributed r.v. with mean $a$, variance $b$.
\begin{lemma}\label{gaussbound00}
For all $x > 0$, 
$$(2 \pi)^{-\frac{1}{2}} \frac{x}{x^2 + 1} \exp( - \frac{x^2}{2} ) < \Phi^c(x) < (2 \pi)^{-\frac{1}{2}} x^{-1} \exp( - \frac{x^2}{2} ).$$
\end{lemma}
\ \\We now complete the proof of Theorem\ \ref{largebtheorem}.
\proof{Proof:}[Proof of Theorem\ \ref{largebtheorem}]
It follows from Corollary\ \ref{mainasymptotic2} (with $\eta = \frac{B}{2}$) and a union bound that for all $B > 0$, $\limsup_{n \rightarrow \infty} \pr\bigg( Q^n_B(\infty) \geq n \bigg)$ is at most
$$
 \sum_{k=0}^{\infty} \pr\bigg( \sup_{\mu^{-1} k \leq t \leq \mu^{-1} (k+1)} \big( {\mathcal Z}(t) - \frac{B}{2} \mu t \big) \ \ \geq\ \ \frac{B}{2} \bigg),
$$
which by a union bound and the stationary increments property of ${\mathcal Z}$ is at most
\begin{eqnarray}
&\ & \pr\bigg( \sup_{0 \leq t \leq \mu^{-1}} {\mathcal Z}(t)\ \ \geq\ \ \frac{B}{2} \bigg) \label{sumunion2}
\\&\ &\ \ +\ \ \sum_{k=1}^{\infty} \pr\bigg( {\mathcal Z}(\mu^{-1} k) \geq \frac{B}{4} (k + 1) \bigg) \label{sumunion3}
\\&\ &\ \ +\ \ \sum_{k=1}^{\infty} \pr\bigg( \sup_{0 \leq t \leq \mu^{-1}} {\mathcal Z}(t) \geq \frac{B}{4} (k + 1) \bigg) \label{sumunion4}
\end{eqnarray}
It follows from Lemmas\ \ref{intervalbound} - \ref{gaussbound00} that for all $\epsilon > 0$, there exists $B_{\epsilon} < \infty$ such that for all $B \geq B_{\epsilon}$, 
\begin{equation}\label{boundb0aa}
(\ref{sumunion2})\ \textrm{is at most}\ \exp\bigg( - (1 - \epsilon) \frac{1}{8} (2 c^2_S + c^2_A + 1)^{-1} B^2 \bigg);
\end{equation}
and
\begin{equation}\label{boundb0aaa}
(\ref{sumunion4})\ \ \ \textrm{is at most}\ \ \ \sum_{k=1}^{\infty}  \exp\bigg( - (1 - \epsilon) \frac{1}{32} (2 c^2_S + c^2_A + 1)^{-1} (k+1)^2 B^2 \bigg).
\end{equation}
Furthermore, it follows from the basic properties of Gaussian r.v.s, and Lemmas\ \ref{boundvarzztop} -\ \ref{gaussbound00}, that 
\begin{equation}\label{boundb0aaaa}
(\ref{sumunion3})\ \ \ \textrm{is at most}\ \ \ \sum_{k=1}^{\infty} \exp\Bigg( - \frac{1}{32} (2 c^2_S + c^2_A + 1)^{-1} (k+1) B^2 \Bigg).
\end{equation}
Combining (\ref{boundb0aa}) - (\ref{boundb0aaaa}), bounding $(k+1)^2$ (from below) by $k+1$ in (\ref{boundb0aaa}), summing all relevant geometric series, and taking the double-limit as $\epsilon \downarrow 0$ and $B \rightarrow \infty$ completes the proof.  $\Halmos$.
\endproof

\section{Proof of bound for s.s.p.d. as $B \rightarrow 0$.}\label{Bzero}
In this section we complete the proof of Theorem\ \ref{smallbtheorem}.  We note that the analysis used to treat the case $B \rightarrow \infty$, i.e. letting $\delta = \infty$ in Theorem\ \ref{mainasymptotic}, will not suffice here.  The underlying reason for this is that when $\delta = \infty$, $\eta$ must be chosen in such a way that stability is maintained, i.e. $\eta \in [0,B)$.  In that case, the bound of Corollary\ \ref{mainasymptotic2} requires a stochastic process (i.e. ${\mathcal Z}$) with drift on the order of $- B$ to be within distance on the order of $B$ of zero.  As $B \rightarrow 0$, such a bound could (at best) demonstrate that the s.s.p.d. is bounded away from unity by a quantity on the order of $B^2$ (e.g. by considering the case of Brownian motion), while the actual bound we will show demonstrates that the s.s.p.d. is bounded away from unity by a much larger quantity (on the order of $B$ as $B \rightarrow 0$).
\\\indent Instead, we will select a non-trivial value for $\delta$.  This will allow us to select a value for $\eta$ in Theorem\ \ref{mainasymptotic} which is independent of $B$, as when $\delta < \infty$, Theorem\ \ref{mainasymptotic} yields non-trivial bounds even when the choice of $\eta$ drives the queue into instability (i.e. $\eta > B$).  Our proof will rely on a careful analysis of the behavior of ${\mathcal Z}(t) - B t$ and its supremum and hitting times, ultimately using the celebrated Slepian's lemma to bound the associated r.v.s by those of simpler Gaussian processes (i.e. the combination of a Brownian motion and an Ornstein-Uhlenbeck process).  We begin by reviewing several results from the theory of Gaussian processes for later use in our analysis. 

\subsection{Slepian's lemma.}\label{slepiansubsec}
Key to our bounds will be the celebrated Slepian's lemma, which will allow us to compare the supremum of ${\mathcal Z}$ to that of a much simpler Gaussian process.  We begin by formally stating a particular variant of Slepian's lemma (see \citet{adler1990introduction}), and a useful corollary.
\begin{lemma}[Slepian's lemma]\label{slepian2a}
Let $T$ denote any finite union of intervals (possibly unbounded) of $\reals^+$.  Let ${\mathcal X}$ and ${\mathcal Y}$ denote any two continuous zero-mean Gaussian processes s.t.
$\E[{\mathcal X}^2(t)] = \E[{\mathcal Y}^2(t)]$ for all $t \in T$, and $\E[{\mathcal X}(s){\mathcal X}(t)] \geq \E[{\mathcal Y}(s){\mathcal Y}(t)]$ for all $s,t \in T$.
Let $g$ denote any function which is continuous on $T$.  Then 
$$\pr\bigg( \sup_{t \in T} \big( {\mathcal X}(t) - g(t) \big) \leq 0 \bigg) \geq \pr\bigg( \sup_{t \in T} \big( {\mathcal Y}(t) - g(t) \big) \leq 0 \bigg).$$
\end{lemma}
\begin{corollary}\label{slepian2}
Let $T^1$ and $T^2$  denote any disjoint intervals of $\reals^+$.  Let ${\mathcal X}$ denote any continuous zero-mean Gaussian processes s.t.
$\E[{\mathcal X}(s){\mathcal X}(t)] \geq 0$ for all $s,t \geq 0$.  Let $g$ denote any function which is continuous on $T^1 \bigcup T^2$.  Then 
$$\pr\bigg( \sup_{t \in T^1 \bigcup T^2} \big( {\mathcal X}(t) - g(t) \big) \leq 0 \bigg) \geq \pr\bigg( \sup_{t \in T^1} \big( {\mathcal X}(t) - g(t) \big) \leq 0 \bigg)
\ \times\ \pr\bigg( \sup_{t \in T^2} \big( {\mathcal X}(t) - g(t) \big) \leq 0 \bigg).$$
\end{corollary}
To apply this corollary, we will need the following result asserting the non-negativity of the covariance of ${\mathcal D}$, whose proof we defer to the appendix (see Subsection\ \ref{proofposcovar0}).
\begin{lemma}\label{poscovar0}
For all $s,t \geq 0$, it holds that $\E[{\mathcal D}(s) {\mathcal D}(t)] \geq 0$.
\end{lemma}
Finally, we state a result closely related to Slepian's inequality, the so-called Sudakov-Fernique Inequality (see \citet{adler1990introduction}), as well as an important corollary bounding the expected value of the supremum of ${\mathcal D}$ (whose proof we defer to the appendix, see Subsection\ \ref{proofboundDsupmean}).
\begin{lemma}[Sudakov-Fernique Inequality]\label{sudakov}
For $T \in \reals^+$, let ${\mathcal X}$ and ${\mathcal Y}$ denote any two continuous zero-mean Gaussian processes s.t.
$\E[\big({\mathcal X}(t) - {\mathcal X}(s)\big)^2] \leq \E[\big({\mathcal Y}(t) - {\mathcal Y}(s)\big)^2]$ for all $s,t \in [0,T]$.  Then
$$\E[\sup_{t \in [0,T]} {\mathcal X}(t) ] \leq \E[\sup_{t \in [0,T]} {\mathcal Y}(t) ].$$
\end{lemma}
\begin{corollary}\label{boundDsupmean}
For $T \in \reals^+$, $\E[\sup_{t \in [0,T]} {\mathcal D}(t) ] \leq \bigg(\frac{2 \mu}{\pi}\big(2 c^2_S + 1 \big) T\bigg)^{\frac{1}{2}}.$
\end{corollary}

\subsection{Properties of Brownian motion and the Ornstein-Uhlenbeck process.}\label{bmou}
In this subsection we review several properties of Brownian motion (B.m.) and the Ornstein-Uhlenbeck (O.U.) process.  In our analysis, we will construct a process whose supremum bounds that of ${\mathcal Z}(t) - B t$ (in an appropriate sense) by taking a weighted combination of independent B.m. and O.U. processes.  For this reason, we will need some relevant results pertaining to the supremum and hitting times of these processes, which we now review.
\subsubsection{Brownian motion.}
For $b > 0$, let $\lbrace {\mathcal B}^b_i(t), i \geq 1 \rbrace$ denote a collection of mutually independent B.m.s initialized to $b$; namely, the continuous Gaussian process s.t. $\E[{\mathcal B}^b(t)] = b$,  $V[{\mathcal B}^b(s),{\mathcal B}^b(t)] = s$ for all $0 \leq s \leq t$.  
As a notational convenience, we will sometimes denote ${\mathcal B}^b_1$ by ${\mathcal B}^b$.  We note that ${\mathcal A}$ is distributed (at the process-level) as $(\mu c^2_A)^{\frac{1}{2}} {\mathcal B}^0$.  We now state several results pertaining to the supremum and hitting times of B.m.  Although some of these results are generally well-known, others are not (especially Lemma\ \ref{brownprops00}.(\ref{brownbig1}) whose proof is non-trivial and involves several intricate properties of conditioned B.m.), and we defer all relevant proofs to the appendix.   For a stochastic process $Z$ and $a \in \bbr$, let $\tau^a_Z$ denote the first hitting time of $Z$ to $a$, where we let $\tau^a_Z = \infty$ if no such time exists.
\begin{lemma}\label{brownprops00}
\ \begin{enumerate}[(i)]
\item For all $t,x > 0$, $\pr\big( \sup_{0 \leq s \leq t} {\mathcal B}^0(s) > x \big) = 2 \Phi^c(x t^{-\frac{1}{2}})$. \label{supform}
\item For all $x,y \in \bbr$ and $z > 0$, $\pr\big( \tau^y_{{\mathcal B}^x} > z\big) = \int_z^{\infty} |y - x| \exp\big( - \frac{(x-y)^2}{2 s} \big) (2 \pi s^3)^{-\frac{1}{2}} ds \leq |y - x| z^{-\frac{1}{2}}$. \label{hitdist00} 
\item For all $c_1,c_2 > 0$, $\pr\big( \tau^{c_1}_{{\mathcal B}^0} < \tau^{- c_2}_{{\mathcal B}^0} \big) = \frac{c_2}{c_1+c_2}.$\label{hitwhich}
\item For all $c,x > 0$, $\pr\bigg(\sup_{t \geq 0}\big( {\mathcal B}^0(t) - c t \big) > x \bigg) = \exp( - 2 c x )$.\label{supdrift1}
\item For all $\epsilon \in (0,1)$, $\gamma > 0$, and $T \geq \big( \gamma^2 \epsilon (1 - \epsilon)^3 \big)^{-1}$, 
$$\pr\bigg(\sup_{0 \leq t \leq T}\big( {\mathcal B}^0(t) + \gamma t \big) \leq \epsilon \gamma T \bigg) \geq (2 \pi)^{-\frac{1}{2}} \frac{\epsilon}{1 - \epsilon^2} (\gamma T^{\frac{1}{2}})^{-1} \exp\big( - \frac{1}{2} (1 - \epsilon)^2 \gamma^2 T \big).$$
\label{supdrift11b}
\item For all $C \geq 1, b \geq 10^{13} C^{10}, \gamma \in \big(0, \frac{10^3 C}{4 b}\big)$, 
$$\pr\Bigg(\sup_{t \geq 0} \bigg( {\mathcal B}^{-b}(t) + C \log^{\frac{1}{2}}(3 + 3 t) - \gamma t\bigg)\ \leq 0 \Bigg) \geq \frac{b \gamma}{10^5 C}.$$
\label{brownbig1}
\end{enumerate}
\end{lemma}
\subsubsection{Ornstein-Uhlenbeck process.}
For any $\rho > 0$, let ${\mathcal U}^{\rho}$ denote the centered, stationary O.U. process whose correlations decay exponentially (over time) at rate $\rho$.  Namely, ${\mathcal U}^{\rho}$ is the continuous Gaussian process s.t. $\E[{\mathcal U}^{\rho}(t)] = 0, V[{\mathcal U}^{\rho}(s), {\mathcal U}^{\rho}(t)] = \exp\big(- \rho(t-s) \big)$ for all $0 \leq s \leq t$.  For a review of the basic properties of O.U processes (e.g. existence, continuity), we refer the interested reader to \citet{Doob.42}.  The main result w.r.t. O.U. processes that we will need in our analysis is a law of the iterated logarithm (l.i.l)-type result.  Although such results are generally well-known (also in considerably greater generality, see \citet{Marcus.72}), we now state a particular variant with associated explicit bounds, whose proof we include in the appendix for completeness.
\begin{lemma}\label{ousup2}
For any fixed $\rho > 0$, 
$$\pr\bigg( \sup_{t \geq 0} \frac{|{\mathcal U}^{\rho}(t)|}{6 \exp(\rho) \log^{\frac{1}{2}}(3 + 3 t)} \leq 1 \bigg) \geq \frac{1}{2}.$$
Namely, with probability at least $\frac{1}{2}$, $|{\mathcal U}^{\rho}(t)|$ is (for all times $t \geq 0$) dominated by $6 \exp(\rho) \log^{\frac{1}{2}}(3 + 3 t)$.
\end{lemma}

\subsection{Critical comparison result.} We now formally state the comparison result which will allow us to compare the supremum of ${\mathcal Z}$ (as well as ${\mathcal D}$) to that of a much simpler Gaussian process.  The result follows from an in-depth analysis of the covariance of ${\mathcal D}$ (and hence that of equilibrium renewal processes), and we defer the proof to the appendix.
Recall that $\alpha_S = \mu^3 \bigg( \E[S^2] + 2 \E[S^3] + \frac{3}{8} \mu \big(\E[S^2]\big)^2 \bigg).$
Suppose $N$ is a $N(0,1)$ r.v. independent of ${\mathcal D}$, and that ${\mathcal B}^0$ and ${\mathcal U}^{1.5}$ are an independent B.m. (initialized at 0) and stationary O.U. process (as defined above).  
\begin{theorem}\label{compareprocess}
There exists a continuous function $f_S: \reals^+ \rightarrow \reals$, depending only on the distribution of $S$, s.t. $\sup_{t \geq 0} |f_S(t)| \leq \alpha_S$, and with the following property.  For $x \geq \alpha_S$, let ${\mathcal D}'_x$ denote the centered Gaussian process $(\mu c^2_S)^{\frac{1}{2}} {\mathcal B}^0 + \big(f_S + 3 x \big)^{\frac{1}{2}} {\mathcal U}^{1.5}$, i.e. for all $t \geq 0$, 
$${\mathcal D}'_x(t) = (\mu c^2_S)^{\frac{1}{2}} {\mathcal B}^0(t) + \big( f_S(t) + 3 x \big)^{\frac{1}{2}} {\mathcal U}^{1.5}(t).$$  
Then the centered Gaussian process ${\mathcal D} + (3 x)^{\frac{1}{2}} N$, i.e. the process which equals ${\mathcal D}(t) + (3 x)^{\frac{1}{2}} N$ for all $t \geq 0$, and the Gaussian process ${\mathcal D}'_x$ satisfy the conditions of Slepian's inequality in such a way that the supremum of ${\mathcal D}'_x$ stochastically dominates that of ${\mathcal D} + (3 x)^{\frac{1}{2}} N$.  More formally, for all $t \geq 0$, $V[{\mathcal D}'_x(t)] = V[{\mathcal D}(t)] + 3 x$ (i.e. $f_S(t) = V[{\mathcal D}(t)] - \mu c^2_S t$); and for all $0 \leq s < t$, 
$$V[{\mathcal D}'_x(s), {\mathcal D}'_x(t)] \leq V[{\mathcal D}(s) + (3 x)^{\frac{1}{2}} N, {\mathcal D}(t) + (3 x)^{\frac{1}{2}} N] = V[{\mathcal D}(s), {\mathcal D}(t)] + 3 x.$$
\end{theorem}

\subsection{Proof of Theorem\ \ref{smallbtheorem}.}
In this subsection, we complete the proof of Theorem\ \ref{smallbtheorem}. 
\proof{Proof:}[Proof of Theorem\ \ref{smallbtheorem}]
Suppose w.l.o.g. that $B \leq 1$.  Let 
$$\delta \stackrel{\Delta}{=} \bigg(200 + \alpha_S + \mu^2 + \mu^{-2} + \sigma_S^2 + \sigma_S^{-2} + \sigma_A^2 + \sigma_A^{-2} + c_A^2 + c_A^{-2} + c_S^2 + c_S^{-2}\bigg)^4.$$
It follows from Theorem\ \ref{mainasymptotic} and taking complements that
\begin{equation}\label{mainasymptoticcor}
\liminf_{n \rightarrow \infty} \pr\bigg( \big( Q^n_B(\infty) - n \big)n^{-\frac{1}{2}} < 0 \bigg) \geq
\pr\bigg( \big({\mathcal E}^{\delta,\delta}_{B,\infty}(0)\big)^c \bigg).
\end{equation}
Note that by construction of $\alpha_S$ and $\delta$, it holds that $\delta > \mu E[S^2] = 2 \E[R(S)]$, and thus $\pr\big( R(S) \leq \delta \big) \geq \frac{1}{2}$ by Markov's inequality.
It then follows from Lemma\ \ref{poscovar0} and Corollary\ \ref{slepian2} that
$\pr\bigg( \big({\mathcal E}^{\delta,\delta}_{B,\infty}(0)\big)^c \bigg)$ is at least
\begin{eqnarray}
&\ &\ \pr\Bigg( \sup_{0 \leq t \leq \delta} \big( {\mathcal Z}(t) + \delta \mu t \big) \ \ \leq\ \frac{1}{2} \delta \Bigg)\label{eqotherdir2}
\\&\indent&\ \ \ \ \times\ \ \ \ \pr\Bigg(\sup_{t \geq \delta} \big( {\mathcal Z}(t) - B \mu t \big)\ \leq - \delta^2 \mu \Bigg).\label{eqotherdir3}
\end{eqnarray}
We first bound (\ref{eqotherdir2}) from below.  By a union bound, the independence of ${\mathcal A}$ and ${\mathcal D}$, and the fact that ${\mathcal A}$ has the same distribution as $(\mu c^2_A)^{\frac{1}{2}} {\mathcal B}^0$, we conclude that (\ref{eqotherdir2}) is at least
\begin{eqnarray}
&\ &\ \pr\bigg( \sup_{0 \leq t \leq \delta} {\mathcal D}(t) \leq \frac{1}{4} \delta \bigg) \label{smallfact3a}
\\&\ &\ \ \ \times\ \ \ \pr\bigg(\sup_{0 \leq t \leq \delta} \big( {\mathcal B}^0(t) + (\mu c^2_A)^{-\frac{1}{2}}\delta \mu t \big) \leq \frac{1}{4} (\mu c^2_A)^{-\frac{1}{2}} \delta \bigg).\label{smallfact3b}
\end{eqnarray}
We now bound (\ref{smallfact3a}) and (\ref{smallfact3b}) from below, beginning with (\ref{smallfact3a}).  By our definition of $\delta$, it is easily verified that $\delta \geq 8 \bigg( \frac{2 \mu}{\pi} (2 c^2_S + 1) \delta \bigg)^{\frac{1}{2}}$.  It thus follows from 
Corollary\ \ref{boundDsupmean} that $\delta \geq 8 \E[\sup_{t \in [0,\delta]} {\mathcal D}(t)]$.  Hence by Markov's inequality,
\begin{equation}\label{smallfact4}
(\ref{smallfact3a})\ \ \ \textrm{is at least}\ \ \ \frac{1}{2}.
\end{equation}
We next bound (\ref{smallfact3b}) from below using Lemma\ \ref{brownprops00}.(\ref{supdrift11b}), i.e. by setting (in the notation of that lemma) $T = \delta, \gamma =  (\mu c^2_A)^{-\frac{1}{2}}\delta \mu, \epsilon = \frac{ \frac{1}{4} (\mu c^2_A)^{-\frac{1}{2}} \delta }{(\mu c^2_A)^{-\frac{1}{2}}\delta^2 \mu} = \frac{1}{4 \mu \delta}$.
The assumptions needed to apply Lemma\ \ref{brownprops00}.(\ref{supdrift11b}) may be easily verified to hold from the definition of $\delta$. We conclude from Lemma\ \ref{brownprops00}.(\ref{supdrift11b}) that 
\begin{equation}\label{smallfact5bb}
(\ref{smallfact3b})\ \ \ \textrm{is at least}\ \ \ (32 \pi)^{-\frac{1}{2}} \mu^{-\frac{3}{2}} c_A \delta^{-\frac{5}{2}} \exp\big( - \frac{1}{2} \mu c_A^{-2} \delta^3 \big).
\end{equation}
As it is easily verified from our definition of $\delta$ that $\frac{1}{2} (32 \pi)^{-\frac{1}{2}} \mu^{-\frac{3}{2}} c_A \geq \delta^{-1}$, and $\frac{1}{2} \mu c_A^{-2} \leq \delta$, we further conclude from (\ref{smallfact5bb}), (\ref{smallfact3a}) - (\ref{smallfact4}), and the exponential inequality that 
\begin{equation}\label{smallfact6}
(\ref{eqotherdir2})\ \ \ \textrm{is at least}\ \ \ \exp\big( - 2 \delta^4 \big).
\end{equation}
We next bound (\ref{eqotherdir3}) from below.  Let $N$ denote a normally distributed r.v. with 0 mean and unit variance, and ${\mathcal B}^0_1$ a B.m, with $N$, ${\mathcal B}^0_1$, ${\mathcal D}$ mutually independent.
Note that (\ref{eqotherdir3}) equals
\begin{eqnarray}
\ &\ &\ \ \ \pr\Bigg(\sup_{t \geq \delta} \big( {\mathcal A}(t) + {\mathcal D}(t) - B \mu t \big)\ \leq - \mu \delta^2 \Bigg)\nonumber
\\ &\ &\ \ \ =\ \ \ \ \pr\Bigg(\sup_{t \geq \delta} \big( {\mathcal A}(\delta) + \big({\mathcal A}(t) - {\mathcal A}(\delta)\big) + {\mathcal D}(t) - B \mu t \big)\ \leq - \mu \delta^2 \Bigg)\nonumber
\\ &\ &\ \ \ =\ \ \ \ \pr\Bigg(\sup_{t \geq \delta} \big( (\mu c^2_A \delta)^{\frac{1}{2}} N + {\mathcal D}(t) + (\mu c^2_A)^{\frac{1}{2}} {\mathcal B}^0_1(t - \delta) - B \mu t \big)\ \leq - \mu \delta^2 \Bigg)\label{smallfact6}.
\end{eqnarray}
Let ${\mathcal U}^{1.5}$ be a stationary O-U process (as defined in Subsection\ \ref{bmou}), and ${\mathcal B}^0_2, {\mathcal B}^0_3, {\mathcal B}^0_4$ be B.m.s, independent of $N, {\mathcal B}^0_1$ (and one-another).  Then it follows from Lemma\ \ref{slepian2a}, Theorem\ \ref{compareprocess}, and the fact that our definition of $\delta$ implies that $\frac{1}{3}\mu c^2_A \delta \geq \alpha_S$, that (\ref{smallfact6}) is at least 
\begin{equation}\label{smallfact7}
\pr\Bigg(\sup_{t \geq \delta} \big( (\mu c^2_S)^{\frac{1}{2}} {\mathcal B}^0_2(t) + \big( f_S(t) + \mu c^2_A \delta \big)^{\frac{1}{2}} {\mathcal U}^{1.5}(t)  + (\mu c^2_A)^{\frac{1}{2}} {\mathcal B}^0_1(t - \delta) - B \mu t \big)\ \leq -  \mu \delta^2 \Bigg).
\end{equation}
Note that by the basic properties of B.m., (\ref{smallfact7}) is unchanged if we replace ${\mathcal B}^0_2(t)$ by $\delta^{\frac{1}{2}} N + {\mathcal B}^0_3(t - \delta)$, and then further replace the resulting term $(\mu c^2_S)^{\frac{1}{2}} {\mathcal B}^0_3(t - \delta) + (\mu c^2_A)^{\frac{1}{2}} {\mathcal B}^0_1(t - \delta)$ by $\big( \mu( c^2_S + c^2_A) \big)^{\frac{1}{2}} {\mathcal B}^0_4(t - \delta)$.  Combining with the stationarity of ${\mathcal U}^{1.5}$ and the fact that $B \mu t \geq B \mu (t - \delta)$, we conclude that (\ref{smallfact6}) is at least
$$\pr\Bigg(\sup_{t \geq \delta} \bigg( (\mu c^2_S \delta)^{\frac{1}{2}} N + \big(\mu (c^2_S + c^2_A)\big)^{\frac{1}{2}} {\mathcal B}^0_4(t - \delta) + \big( f_S(t) + \mu c^2_A \delta \big)^{\frac{1}{2}} {\mathcal U}^{1.5}(t - \delta) - B \mu (t - \delta)\big)\ \leq - \mu \delta^2 \Bigg),
$$ 
which itself equals
\begin{equation}\label{smallfact8}
\pr\Bigg(\sup_{t \geq 0} \bigg( (\mu c^2_S \delta)^{\frac{1}{2}} N + \big(\mu (c^2_S + c^2_A)\big)^{\frac{1}{2}} {\mathcal B}^0_4(t) + \big( f_S(t + \delta) + \mu c^2_A \delta \big)^{\frac{1}{2}} {\mathcal U}^{1.5}(t) - B \mu t\big)\ \leq - \mu \delta^2 \Bigg).
\end{equation}
Let $C \stackrel{\Delta}{=}50 (\frac{\delta \sigma^2_A}{\sigma^2_S + \sigma^2_A})^{\frac{1}{2}}$, and $\gamma_B \stackrel{\Delta}{=} \big(\mu(\sigma^2_S + \sigma^2_A)\big)^{-\frac{1}{2}} B$.  By conditioning on the event $\lbrace \sup_{t \geq 0} \frac{|{\mathcal U}^{1.5}(t)|}{6 \exp(1.5) \log^{\frac{1}{2}}(3 + 3 t)} \leq 1 \rbrace$, applying Lemma\ \ref{ousup2} (to bound the corresponding probability) and Theorem\ \ref{compareprocess} (which asserts that $\sup_{t \geq 0} |f_S(t)| \leq \alpha_S \leq \mu c^2_A \delta $), and simplying some straightforward algebra, we conclude that (\ref{smallfact8}), and hence (\ref{eqotherdir3}), is at least
\begin{equation}\label{smallfact9}
\frac{1}{2} \pr\Bigg(\sup_{t \geq 0} \bigg( (\frac{\sigma^2_S \delta}{\sigma^2_S + \sigma^2_A})^{\frac{1}{2}} N + {\mathcal B}^0_4(t) + C \log^{\frac{1}{2}}(3 + 3 t) 
- \gamma_B t \bigg)\ \leq - \big(\mu(\sigma^2_S + \sigma^2_A)\big)^{-\frac{1}{2}} \delta^2 \Bigg).
\end{equation}
Let $b \stackrel{\Delta}{=} 10^{13} \times C^{10}$.  Further conditioning on the event 
$$\bigg\lbrace  (\frac{\sigma^2_S \delta}{\sigma^2_S + \sigma^2_A})^{\frac{1}{2}} N \leq - \big(\mu(\sigma^2_S + \sigma^2_A)\big)^{-\frac{1}{2}} \delta^2 - b \bigg\rbrace,$$
and recalling that ${\mathcal B}^{-b}$ denotes a B.m. initialized at $-b$,
we find that (\ref{eqotherdir3}) is at least
\begin{eqnarray}
\ &\ &\ \ \ \frac{1}{2} \Phi^c\Bigg( (\mu \sigma^2_S)^{-\frac{1}{2}} \delta^{\frac{3}{2}} + b \big(\frac{\sigma^2_S + \sigma^2_A}{\sigma^2_S \delta}\big)^{\frac{1}{2}}\Bigg) \nonumber
\\&\ &\ \ \ \ \ \ \times\ \ \ \pr\Bigg(\sup_{t \geq 0} \bigg( {\mathcal B}^{-b}(t) + C \log^{\frac{1}{2}}(3 + 3 t) - \gamma_B t\bigg)\ \leq 0 \Bigg). \label{lastbig1}
\end{eqnarray}
It is easily verified from our definitions of $\delta, C, b$ that for all sufficiently small $B$, we may apply Lemma\ \ref{brownprops00}.(\ref{brownbig1}) to bound (\ref{lastbig1}).  In particular, there exists $B_0 > 0$ s.t. for all $B \in (0,B_0)$, 
$$
\pr\Bigg(\sup_{t \geq 0} \bigg( {\mathcal B}^{-b}(t) + C \log^{\frac{1}{2}}(3 + 3 t) - \gamma_B t\bigg)\ \leq 0 \Bigg) \geq \frac{b \gamma_B}{10^5 C}.$$
Combining with (\ref{lastbig1}), we conclude that for all $B \in (0,B_0)$, (\ref{eqotherdir3}) is at least
\begin{equation}\label{lastbig2}
\frac{1}{2} \times \Phi^c\Bigg( (\mu \sigma^2_S)^{-\frac{1}{2}} \delta^{\frac{3}{2}} + b \big(\frac{\sigma^2_S + \sigma^2_A}{\sigma^2_S \delta}\big)^{\frac{1}{2}}\Bigg) \times \frac{b \gamma_B}{10^5 C}.
\end{equation}
As our definitions of $\delta, C$, and $b$ imply that $\frac{\sigma^2_S + \sigma^2_A}{\sigma^2_S \delta} \leq 1$ and $(\mu \sigma^2_S)^{-\frac{1}{2}} \delta^{\frac{3}{2}} \leq b$, it further holds that (\ref{eqotherdir3}) is at least
\begin{equation}\label{lastbig2b}
\frac{1}{2} \times \Phi^c(2 b) \times \frac{b \gamma_B}{10^5 C}.
\end{equation}
Combining with (\ref{smallfact6}) and (\ref{mainasymptoticcor}), we conclude that for all $B \in (0,B_0)$, it holds that $\liminf_{n \rightarrow \infty} \pr\bigg( \big( Q^n_B(\infty) - n \big)n^{-\frac{1}{2}} < 0 \bigg)$ is at least
\begin{equation}\label{lastbig3}
\frac{1}{2} \times \Phi^c(2 b) \times \exp\big( - 2 \delta^4 \big) \times \frac{b \gamma_B}{10^5 C}.
\end{equation}
As Lemma\ \ref{gaussbound00} and our definitions of $\delta,C$, and $b$ imply that $\Phi^c(2 b) \geq \frac{1}{6} b^{-1} \exp( - 2 b^2)$, $\frac{\gamma_B}{10^5 C} \geq \exp(- b^2) B$, and $\frac{1}{12} \exp\big( - 2 \delta^4 \big) \geq \exp(- b^2)$, we conclude that 
$\liminf_{n \rightarrow \infty} \pr\bigg( \big( Q^n_B(\infty) - n \big)n^{-\frac{1}{2}} < 0 \bigg)$ is at least $\exp(- 4 b^2) \times B$\ for all\ $B \in (0,B_0)$.  As our definitions of $\delta, C$, and $b$ imply that $C \leq \delta$ and thus $b \leq 10^{13} \delta^{10}$, combining the above completes the proof.
$\Halmos$.
\endproof
\section{Proof of large deviations result for number of idle servers.}\label{negx}
In this section we complete the proof of Theorem\ \ref{negxtheorem}.  As in the proof of Theorem\ \ref{smallbtheorem}, here we will have to carefully choose non-trivial values for both $\delta$ and $\eta$ in Theorem\ \ref{mainasymptotic}, since e.g. choosing $\delta = \infty$ would result in instability and trivial bounds.  We note that here, $\delta$ will be selected as some non-trivial value not depending on $x$, while $\eta$ will be selected as an appropriate function of $x$.
\proof{Proof:}[Proof of Theorem\ \ref{negxtheorem}]
Let us fix some $B > 0$, $x < -1$, and let $\alpha_x \stackrel{\Delta}{=} 2 |x| + 6 |x|^{\frac{1}{2}}$.  It follows from Theorem\ \ref{mainasymptotic} that 
$\liminf_{n \rightarrow \infty} \pr\bigg( \big( Q^n_B(\infty) - n \big)n^{-\frac{1}{2}} < x \bigg)$ is at least
$\pr\big({\mathcal E}^{\mu \E[S^2], \alpha_x }_{B,\infty}(x)\big)^c$, which (since $\mu \E[S^2] \geq 2 \E[R(S)]$ and thus 
$x + \alpha_x \pr\big( R(S) \leq \mu \E[S^2] \big) \geq 3 |x|^{\frac{1}{2}}$) is at least 
\begin{equation}\label{negxeq1}
\pr\Bigg( \max \bigg( \sup_{0 \leq t \leq \mu \E[S^2]} \bigg( {\mathcal Z}(t) + \big( \alpha_x - B \big) \mu t \bigg)\ ,\ 
\sup_{t \geq \mu \E[S^2]} \big( {\mathcal Z}(t) - B \mu t \big)\ +\ \alpha_x \mu^2 \E[S^2] \bigg)
\leq 3 |x|^{\frac{1}{2}} \Bigg).
\end{equation}
It then follows from a union bound and the independence of ${\mathcal A}$ and ${\mathcal D}$ that
(\ref{negxeq1}) is at least
\begin{eqnarray}
\ &\ &\ \ \pr\big( \sup_{t \geq 0} ( {\mathcal D}(t) - \frac{1}{2} B \mu t ) \leq |x|^{\frac{1}{2}} \big) \label{negxeq2}
\\&\ &\ \ \times\ \pr\Bigg( \bigg\lbrace \sup_{0 \leq t \leq \mu \E[S^2]} \big( {\mathcal A}(t) + \alpha_x \mu t \big) \leq 2 |x|^{\frac{1}{2}} \bigg\rbrace\ \ \ , \label{negxeq3}
\\&\ &\ \ \ \ \ \ \ \ \ \ \ \ \ \ \ \ \bigg\lbrace \sup_{t \geq \mu \E[S^2]} \big( {\mathcal A}(t) - \frac{1}{2} B \mu t \big)
\leq - \big( \mu^2 \E[S^2] \alpha_x \ -\ 2 |x|^{\frac{1}{2}} \big) \bigg\rbrace \Bigg). \nonumber
\end{eqnarray}
We now bound (\ref{negxeq3}).  Let $\gamma_{B} \stackrel{\Delta}{=} \frac{1}{2} (\mu c^2_A)^{-\frac{1}{2}} B \mu$.
Note that by the independent increments property of B.m., we may construct ${\mathcal B}^0_2$ and ${\mathcal A}$  
as independent processes on a common probability space s.t. for all $t \geq \mu \E[S^2]$, ${\mathcal A}(t) = 
{\mathcal A}(\mu \E[S^2]) + (\mu c^2_A)^{\frac{1}{2}} {\mathcal B}^0_2(t - \mu \E[S^2])$.  It follows that (\ref{negxeq3}) equals
\begin{eqnarray*}
\ &\ &\ \pr\Bigg( \bigg\lbrace \sup_{0 \leq t \leq \mu \E[S^2]} \big( {\mathcal A}(t) + \alpha_x \mu t \big) \leq 2 |x|^{\frac{1}{2}} \bigg\rbrace\ \ \ ,
\\&\ &\ \ \ \ \ \ \ \ \bigg\lbrace (\mu c^2_A)^{-\frac{1}{2}} {\mathcal A}(\mu \E[S^2]) + \sup_{t \geq 0} \big( {\mathcal B}^0_2(t) - 
\gamma_{B} (t + \mu \E[S^2]) \big)
\leq 
- (\mu c^2_A)^{-\frac{1}{2}} \big( \mu^2 \E[S^2] \alpha_x\ -\ 2 |x|^{\frac{1}{2}} \big) \bigg\rbrace \Bigg),
\end{eqnarray*}
which is itself at least
\begin{equation}\label{negxeq4}
\pr\Bigg( \sup_{0 \leq t \leq \mu \E[S^2]} \big( {\mathcal A}(t) + \alpha_x \mu t \big) \leq |x|^{\frac{1}{2}}\ \ \ ,\ \ \ \sup_{t \geq 0} \big( {\mathcal B}^0_2(t) - \gamma_{B} t \big)
\leq (\mu c^2_A)^{-\frac{1}{2}}|x|^{\frac{1}{2}} \Bigg), 
\end{equation}
the final inequality following from the fact that $\sup_{0 \leq t \leq \mu \E[S^2]} \big( {\mathcal A}(t) + \alpha_x \mu t \big) \leq |x|^{\frac{1}{2}}$ implies ${\mathcal A}(\mu \E[S^2]) \leq - \big( \mu^2 \E[S^2] \alpha_x - |x|^{\frac{1}{2}} \big)$.
Combining with the independence of ${\mathcal A}$ and ${\mathcal B}^0_2$, and the fact that ${\mathcal A}$ is distributed as $(\mu c^2_A)^{\frac{1}{2}} {\mathcal B}^0$, we conclude that (\ref{negxeq4}) is at least
\begin{eqnarray}
\ &\ &\ \pr\bigg( \sup_{0 \leq t \leq \mu \E[S^2]} \big( {\mathcal B}^0(t) + 
(\mu c^2_A)^{-\frac{1}{2}} \alpha_x \mu t \big) \leq (\mu c^2_A)^{-\frac{1}{2}} |x|^{\frac{1}{2}} \bigg) \label{negxeq5}
\\&\ &\ \ \ \times\ \ \ \pr\bigg( \sup_{t \geq 0} \big( {\mathcal B}^0(t) - \frac{1}{2} (\mu c^2_A)^{-\frac{1}{2}} B \mu t \big)
\leq  (\mu c^2_A)^{-\frac{1}{2}}|x|^{\frac{1}{2}} \bigg). \label{negxeq6}
\end{eqnarray}
We next bound (\ref{negxeq5}) using Lemma\ \ref{brownprops00}.(\ref{supdrift11b}), by setting (in the notation of that lemma) $T = \mu \E[S^2], \gamma = (\mu c^2_A)^{-\frac{1}{2}} (2 |x| + 6 |x|^{\frac{1}{2}}) \mu, \epsilon = \frac{ (\mu c^2_A)^{-\frac{1}{2}} |x|^{\frac{1}{2}}}{ (\mu c^2_A)^{-\frac{1}{2}} (2 |x| + 6 |x|^{\frac{1}{2}}) \mu^2 \E[S^2]} = \big( (2|x|^{\frac{1}{2}} + 6) \mu^2 \E[S^2]\big)^{-1}$.  It is easily verified that there exists $x_0 > - \infty$ s.t. $x < x_0$ implies that the assumptions of Lemma\ \ref{brownprops00}.(\ref{supdrift11b}) are met, and (\ref{negxeq5}) is at least
\begin{equation}\label{negxeq7}
(2 \pi)^{-\frac{1}{2}} \big( (2|x|^{\frac{1}{2}} + 6) \mu^2 \E[S^2]\big)^{-1} \big( \sigma^{-1}_A (2 |x| + 6 |x|^{\frac{1}{2}})(\E[S^2])^{\frac{1}{2}} \big)^{-1}
\exp\big( - \frac{1}{2} \E[S^2] \sigma^{-2}_A  (2 |x| + 6 |x|^{\frac{1}{2}})^2 \big).
\end{equation}
As it follows from Theorem\ \ref{oldandweak}.(\ref{varlimitplus}) that 
$$
\lim_{x \rightarrow -\infty} \pr\big( \sup_{t \geq 0} ( {\mathcal D}(t) - \frac{1}{2} B \mu t ) \leq |x|^{\frac{1}{2}} \big) = 1,$$
and
$$
\lim_{x \rightarrow -\infty} \pr\bigg( \sup_{t \geq 0} \big( {\mathcal B}^0(t) - \frac{1}{2} (\mu c^2_A)^{-\frac{1}{2}} B \mu t \big) \leq  (\mu c^2_A)^{-\frac{1}{2}}|x|^{\frac{1}{2}} \bigg) = 1,$$
combining (\ref{negxeq1}) - (\ref{negxeq7}) and taking limits in (\ref{negxeq7}) completes the proof.  $\Halmos$
\endproof
\section{Comparison to other bounds from the literature.}\label{4comparesec}
In this section we compare our results to known results for the $GI/M/n$, $GI/D/n$, and $M/GI/\infty$ queues, showing that our main results are tight, in an appropriate sense.  
\subsection{$GI/M/n$ queue.}
Let us define $\alpha(x) \stackrel{\Delta}{=}
\big( 1 + x \Phi(x) \phi^{-1}(x) \big)^{-1}$.  In \citet{Whitt.05d}, the authors prove the following, which generalizes (and corrects) the results for the $GI/M/n$ queue given in \citet{HW.81}.
We restrict our discussion to the $GI/M/n$ setting (as opposed to the more general $GI/H^*_2/n$ setting treated in \citet{Whitt.05d}) as the general $GI/H^*_2/n$ setting does not satisfy our $T_0$ assumptions, as in that case the service distribution puts strictly positive probability at the origin.
\begin{lemma}[\citet{Whitt.05d}]\label{hwtheorem11}
Suppose that $\lbrace {\mathcal Q}^n_B, n \geq 1 \rbrace$ is a sequence of $GI/M/n$ queues satisfying the H-W and $T_0$ assumptions.  Also suppose that $\E[A^3] < \infty$.  Let $z \stackrel{\Delta}{=} \frac{1}{2} (c^2_A + 1)$.  Then for all $x \geq 0$, 
$$\lim_{n \rightarrow \infty} \pr\bigg( \big(Q^n_B(\infty) - n\big)n^{-\frac{1}{2}} \geq x \bigg) = 
\alpha(B z^{-\frac{1}{2}}) \exp( - B z^{-1} x ),$$
and
$$\lim_{n \rightarrow \infty} \pr\bigg( \big(Q^n_B(\infty) - n\big)n^{-\frac{1}{2}} \leq - x \bigg) = 
\big(1 - \alpha(B z^{-\frac{1}{2}})\big) \frac{ \Phi\big( (B-x)z^{-\frac{1}{2}}\big)}{\Phi( B z^{-\frac{1}{2}})}.$$
\end{lemma}
Lemma\ \ref{hwtheorem11} and a straightforward asymptotic analysis (the details of which we omit) then yield the following.  
\begin{corollary} \label{hwcorr1}
Under the same assumptions as Lemma\ \ref{hwtheorem11},
\begin{enumerate}[(i)]
\item $\lim_{B \rightarrow \infty} B^{-2} \log \bigg( \lim_{n \rightarrow \infty} \pr\big( Q^n(\infty) \geq n \big) \bigg) = - \big(c^2_A+1\big)^{-1},$
\item $\lim_{B \rightarrow 0} B^{-1}  \lim_{n \rightarrow \infty} \pr\big( Q^n(\infty) < n \big) = \pi^{\frac{1}{2}} \big(c^2_A+1\big)^{-\frac{1}{2}},$
\item $\lim_{x \rightarrow \infty} x^{-2} \log \Bigg( \lim_{n \rightarrow \infty} \pr\bigg( \big( Q^n(\infty) - n \big)n^{-\frac{1}{2}} \leq - x \bigg) \Bigg) = - \big(c^2_A+1\big)^{-1}.$
\end{enumerate}
\end{corollary}
Alternatively, our main results imply the following bounds.
\begin{observation} \label{hwcorr111}
.\ \ Under the same assumptions as Lemma\ \ref{hwtheorem11}, our main results imply the following bounds for the $GI/M/n$ queue:
\begin{enumerate}[(i)]
\item $\limsup_{B \rightarrow \infty} B^{-2} \log \bigg( \limsup_{n \rightarrow \infty} \pr\big( Q^n(\infty) \geq n \big) \bigg) \leq -\frac{1}{16} \big(c^2_A+3\big)^{-1},$
\item $\liminf_{B \rightarrow 0} B^{-1}  \liminf_{n \rightarrow \infty} \pr\big( Q^n(\infty) < n \big) > 0,$
\item $\liminf_{x \rightarrow \infty} x^{-2} \log \Bigg( \liminf_{n \rightarrow \infty} \pr\bigg( \big( Q^n(\infty) - n \big)n^{-\frac{1}{2}} \leq - x \bigg) \Bigg) \geq - 4 c_A^{-2}.$
\end{enumerate}
\end{observation}
Note that our bounds for the s.s.p.d. (as $B \rightarrow \infty$) and for the large deviations exponent for the number of idle servers are quite close to the true exponents, have a similar functional form, and in some cases even have a uniformly bounded error.  For example, 
$\sup_{c_A \geq 0} \frac{ \big(c^2_A+1\big)^{-1} }{ \frac{1}{16} \big(c^2_A+3\big)^{-1} } = 48$, i.e. our bound for the s.s.p.d. (as $B \rightarrow \infty$) is within a multiplicative factor of 48 of the true value for any value of $c_A$.  For the case $B \rightarrow 0$, we note that our results in principle give a quantitative estimate (see Theorem\ \ref{smallbtheorem}), but here (for clarity of exposition) we have only stated the implication that the corresponding limit is strictly positive, which agrees qualitatively with the true limit.
\subsection{$GI/D/n$ queue.}
In \citet{JMM.04}, the authors prove the following.
\begin{lemma}[\citet{JMM.04}]\label{dettheorem1}
Suppose that $\lbrace {\mathcal Q}^n_B, n \geq 1 \rbrace$ is a sequence of $GI/D/n$ queues satisfying the H-W and $T_0$ assumptions, i.e. the processing times are distributed as some strictly positive constant.  Let $\lbrace X_i, i \geq 1 \rbrace$ denote a set of i.i.d. normally distributed r.v.s with mean $-B$ and variance $c^2_A > 0$.  Let $S_i \stackrel{\Delta}{=} \sum_{k=1}^i X_k , i \geq 0$.  Then for all $x \in \reals$,
$$\lim_{n \rightarrow \infty} \pr\bigg( \big(Q^n_B(\infty) - n\big)n^{-\frac{1}{2}} \geq x \bigg) = 
\pr\big( \sup_{i \geq 1} S_i \geq x \big).$$
\end{lemma}
Lemma\ \ref{dettheorem1} implies the following, whose proof we defer to the appendix.
\begin{corollary}\label{dettheoremcor}
Under the same assumptions as Lemma\ \ref{dettheorem1}:
\begin{enumerate}[(i)]
\item $\lim_{B \rightarrow \infty} B^{-2} \log \bigg( \lim_{n \rightarrow \infty} \pr\big( Q^n(\infty) \geq n \big) \bigg) = - (2 c^2_A)^{-1},$ \label{detcor2a}
\item $\lim_{B \rightarrow 0} B^{-1}  \lim_{n \rightarrow \infty} \pr\big( Q^n(\infty) < n \big) = 2^{\frac{1}{2}} c^{-1}_A,$ \label{detcor2b}
\item $\lim_{x \rightarrow \infty} x^{-2} \log \Bigg( \lim_{n \rightarrow \infty} \pr\bigg( \big( Q^n(\infty) - n \big)n^{-\frac{1}{2}} \leq - x \bigg) \Bigg) = - (2 c^2_A)^{-1}.$\label{detcor2c}
\end{enumerate}
\end{corollary}
Alternatively, our main results imply the following bounds.
\begin{observation} \label{hwcorr222}
.\ \ Under the same assumptions as Lemma\ \ref{dettheorem1}, our main results imply the following bounds for the $GI/D/n$ queue:
\begin{enumerate}[(i)]
\item $\limsup_{B \rightarrow \infty} B^{-2} \log \bigg( \limsup_{n \rightarrow \infty} \pr\big( Q^n(\infty) \geq n \big) \bigg) \leq -\frac{1}{16} \big(c^2_A+1\big)^{-1},$
\item $\liminf_{B \rightarrow 0} B^{-1}  \liminf_{n \rightarrow \infty} \pr\big( Q^n(\infty) < n \big) > 0,$
\item $\liminf_{x \rightarrow \infty} x^{-2} \log \Bigg( \liminf_{n \rightarrow \infty} \pr\bigg( \big( Q^n(\infty) - n \big)n^{-\frac{1}{2}} \leq - x \bigg) \Bigg) \geq - 2 c_A^{-2}.$
\end{enumerate}
\end{observation}
Once again, our bounds are quite close to the true exponents, and have a similar functional form.  Note that here, for any value of $c_A$ our bound for the large deviations exponent for the number of idle servers is always exactly a multiplicative factor of 4 away from the true value.  Also, once again our results for the setting $B \rightarrow 0$ correctly capture the qualitative scaling of the true limit.
\subsection{$M/GI/\infty$ lower bound.}
Suppose that $Q^n_B$ is an $M/GI/n$ queue satisfying the H-W and $T_0$ assumptions.  Let $Z_{n,B}$ denote a Poisson r.v. with mean $\lambda_{n,B}$.  Then it follows from a straightforward infinite-server lower bound, and the well-known properties of the steady-state infinite server queue (see \citet{Tak.62}), that for all $x \in \reals^+$, $\pr\big( Q^n_B(\infty) \leq n - x n^{\frac{1}{2}} \big) \leq \pr\big( Z_{n,B} \leq n - x n^{\frac{1}{2}} \big)$.  
It follows from the Central Limit Theorem that for all $x \in \reals^+$,
$\lim_{n \rightarrow \infty} \pr ( Z_{n,B} \leq n - x n^{\frac{1}{2}} ) = \Phi(B - x)$, and we conclude the following.
\begin{lemma}\label{infinitebound}
Suppose that $\lbrace Q^n_B, n \geq 1 \rbrace$ is a sequence of $M/GI/n$ queues satisfying the H-W and $T_0$ assumptions.  Then
$$\liminf_{B \rightarrow \infty} B^{-2} \log \bigg( \liminf_{n \rightarrow \infty} \pr \big( Q^n_B(\infty) \geq n \big) \bigg) \geq - \frac{1}{2},$$
and
$$
\limsup_{x \rightarrow \infty} x^{-2} \log \Bigg( \limsup_{n \rightarrow \infty} \pr\bigg( \big( Q^n(\infty) - n \big)n^{-\frac{1}{2}} \leq - x \bigg) \Bigg) \leq - \frac{1}{2}.
$$
\end{lemma}
Alternatively, our main results imply the following bounds.
\begin{observation} \label{hwcorr222}
.\ \ Under the same assumptions as Lemma\ \ref{infinitebound}, our main results imply the following bounds for the $M/G/n$ queue:
\begin{enumerate}[(i)]
\item $\limsup_{B \rightarrow \infty} B^{-2} \log \bigg( \limsup_{n \rightarrow \infty} \pr\big( Q^n(\infty) \geq n \big) \bigg) \leq -\frac{1}{32} \big(c^2_S+1\big)^{-1},$
\item $\liminf_{x \rightarrow \infty} x^{-2} \log \Bigg( \liminf_{n \rightarrow \infty} \pr\bigg( \big( Q^n(\infty) - n \big)n^{-\frac{1}{2}} \leq - x \bigg) \Bigg) \geq - 2 \mu^2 \E[S^2].$
\end{enumerate}
\end{observation}
Once again, our bounds are quite close to the true exponents, at least for moderate values of $c^2_S$.  Interestingly, the infinite server lower-bound does not seem to yield any information about the tightness of Theorem\ \ref{smallbtheorem}, since $\Phi(B) \geq \frac{1}{2}$ for all $B \geq 0$.
\section{Conclusion and open questions.}\label{4concsec}
In this paper, we studied the FCFS $GI/GI/n$ queue in the H-W regime, deriving bounds on the s.s.p.d. and number of idle servers.  We proved that under quite general assumptions on the inter-arrival and processing time distributions, e.g. finite third moment, there exists some strictly positive $\epsilon_1,\epsilon_2,\epsilon_3$ (depending on the first three moments of the inter-arrival and processing time distributions) s.t. in the H-W regime, the s.s.p.d. is bounded from above by $\exp\big(-\epsilon_1 B^2\big)$ as the associated excess parameter $B \rightarrow \infty$; and by $1 - \epsilon_2 B$ as $B \rightarrow 0$.  We also prove that the probability of there being more than $x n^{\frac{1}{2}}$ idle servers (in steady-state, for large $n$) is bounded from below by $\exp\big( - \epsilon_3 x^2 \big)$ as $x \rightarrow \infty$.  Furthermore, we used known results for the special cases of Markovian and deterministic processing times, as well as the $M/GI/\infty$ queue, to prove that our bounds are tight, in an appropriate sense.  Our main proof technique was the derivation of new stochastic comparison bounds for the FCFS $GI/GI/n$ queue, which are of a structural nature, and significantly extend the recent work of \citet{GG.10c}, combined with several bounding arguments for Gaussian processes (using e.g. Slepian's inequality).   Our results do not follow from simple comparison arguments to e.g. infinite-server systems or loss models, which would in all cases provide bounds in the opposite direction.
\\\indent This work leaves many interesting directions for future research.  The explicit limits for the s.s.p.d. as $B \rightarrow 0$ and $B \rightarrow \infty$ for Markovian and deterministic processing times are suggestive of some fascinating patterns, and it is an intriguing open problem to more precisely understand the nature of these limits for general distributions.   In addition to (or combined with) the stochastic comparison approach, a direct analysis of these quantities for the true limiting processes identified in e.g. \citet{GM.08}, \citet{PR.10}, \citet{dieker2013positive}, and \citet{aghajani2015ergodicity} remains a potentially fruitful direction for research.  It is also interesting to note that for Markovian and deterministic processing times, $\lim_{B \rightarrow \infty} B^{-2} \log(p_B)$ is actually equal to $\lim_{x \rightarrow \infty} x^{-2} \log\big(F_{\textrm{idle},B}(x)\big)$ - to what extent such a relationship may hold in general is unknown.  Furthermore, these limits may fit into the broader context of insensitivity results in queueing, in which certain quantities depend on the relevant distributions only through limited information (e.g. first two moments), and it is open to further investigate such connections.
\\\indent We believe that the stochastic comparison techniques developed in this paper, and the original work of \citet{GG.10c}, have the potential to shed insight on many other queueing models.  Indeed, one may view our methods as a step towards developing a calculus of stochastic-comparison type bounds for parallel server queues, in which one derives bounds by composing structural modifications (e.g. adding jobs, removing servers, adding servers, etc.) over time.  Three particularly interesting systems to which one might try to apply these methods are queues with abandonments, queues with heavy-tailed processing times, and networks of queues, which arise in various applied settings.  We note that the setting of abandonments and networks are particularly interesting from a stochastic comparison standpoint, as these systems may exhibit certain non-monotonicities (e.g. adding a job may cause other jobs to leave the system sooner), and developing a better understanding of when such systems can be compared will likely require the creation of new tools.  
\\\indent On a related note, understanding the fundamental power of these stochastic comparison techniques, e.g. whether they can yield tractable approximations sufficiently tight to capture the true limits and exponents (as was the case in \citet{GG.10c}), remains an open question.  We note that even in the present paper, it may be that certain of our arguments could be made to yield tighter bounds if one were to more carefully optimize over the given free parameters and more precisely bound various probabilities related to Gaussian processes.  Indeed, at times we sided with simplicity as opposed to pushing the methodologies to their tightest possible results, and a more careful analysis may be a good first step towards understanding the ultimate power of these techniques.  For example, if one carefully optimizes over all relevant parameters, how tight a bound does Theorem\ \ref{mainasymptotic} yield for the s.s.p.d. in general (i.e. without letting $B$ go to $0$ or $\infty$)?  Can stochastic comparison techniques be used to develop significantly tighter bounds for the s.s.p.d. as $B \downarrow 0$, where we note that developing good lower bounds (in addition to better upper bounds) remains an open challenge.  It is also interesting to study the extent to which limiting results such as Theorems\ \ref{largebtheorem}\ -\ \ref{negxtheorem} apply for any fixed $n$, e.g. through a careful pre-limit analysis of Theorem\ \ref{4ubound1}, and refer the reader to \citet{GS.12d}, \citet{braverman2015stein}, and \citet{braverman2016high} for some results along these lines.
\\\indent On a final note, the different qualitative behavior we observe w.r.t. the s.s.p.d. as $B \rightarrow 0$ and $B \rightarrow \infty$ fits into the broader theme of the H-W regime as a ``transition" between a system which behaves like an infinite-server queue ($B \rightarrow \infty$) and a single-server queue ($B \rightarrow 0$).  In the H-W regime, this transition was previously formalized for the case of Markovian processing times by \citet{GG.10b} and \citet{van2011transient}, by proving that there exists $B^* \approx 1.85722$ s.t. the spectral measure of the underlying Markov chain had no jumps (like the single-server queue) for all $B \in (0,B^*)$, and had at least one jump (like the infinite-server queue) for all $B > B^*$.  Furthermore, in the work of \citet{van2011transient}, a certain notion of ``convergence" of the spectral measure to that of the infinite-server queue, as $B \rightarrow \infty$, is demonstrated.  It is an interesting open question to more generally formalize these ``limits within limits" of the H-W regime.
\section{Appendix.}\label{4appsec}
\subsection{Proof of Corollary\ \ref{mainasymptotic2}.}
\proof{Proof:}[Proof of Corollary\ \ref{mainasymptotic2}]
It follows from Theorem\ \ref{mainasymptotic}, combined with the monotonicity of the supremum operator and a union bound, that 
for all $B > 0$, $x \in \reals$, $\eta \in [0,B)$, and $\delta \geq 0$, $\limsup_{n \rightarrow \infty} \pr\bigg( \big( Q^n_B(\infty) - n \big)n^{-\frac{1}{2}} > x \bigg)$ is at most 
\begin{eqnarray*}
&\ &\ \ \ \pr\bigg(  \sup_{t \geq 0} \big( {\mathcal Z}(t) - (B - \eta) \mu t \big)\ \ \geq\ \ x + \eta \pr\big( R(S) \leq \delta \big) \bigg)
\\&\ &\ \ \ \ +\ \ \ \pr\bigg( \sup_{t \geq \delta} \big( {\mathcal Z}(t) - B \mu t \big)\ +\ \eta \mu \delta\ \ \geq\ \sup_{0 \leq t \leq \delta} \big( {\mathcal Z}(t) + (\eta - B) \mu t \big) \bigg),
\end{eqnarray*}
which is itself bounded by
\begin{eqnarray}
&\ &\ \ \ \pr\bigg(  \sup_{t \geq 0} \big( {\mathcal Z}(t) - (B - \eta) \mu t \big)\ \ \geq\ \ x + \eta \pr\big( R(S) \leq \delta \big) \bigg)\label{deltalimit0}
\\&\ &\ \ \ \ \ +\ \ \ \pr\bigg(  \sup_{t \geq \delta} \big( {\mathcal Z}(t) - B \mu t \big)\ \ \geq\ - \eta \mu \delta \bigg).\label{deltalimit2}
\end{eqnarray}
By the stationary increments property of ${\mathcal Z}$ and a union bound, (\ref{deltalimit2}) is at most
$$\pr\bigg( {\mathcal Z}(\delta) - B \mu \delta \geq - \frac{1}{2}(B + \eta) \mu \delta \bigg) + \pr\bigg( \sup_{t \geq 0} \big( {\mathcal Z}(t) - B \mu t \big) \geq \frac{1}{2}(B - \eta)\mu \delta \bigg),$$
which is itself equivalent to 
\begin{equation}\label{deltabound1}
\pr\big( {\mathcal Z}(\delta) \geq \frac{1}{2}(B - \eta)\mu \delta \big) + \pr\bigg( \sup_{t \geq 0} \big( {\mathcal Z}(t) - B \mu t \big) \geq \frac{1}{2}(B - \eta)\mu \delta \bigg).
\end{equation}
It follows from Theorem\ \ref{oldandweak}.(\ref{varlimitplus}) that for any $\epsilon > 0$, there exists $\delta_{\epsilon} < \infty$ s.t. for all $\delta \geq \delta_{\epsilon}$, (\ref{deltabound1}) is at most $\epsilon$.  The corollary then follows from (\ref{deltalimit0}) and the continuity of probability measures, since 
$\lim_{\delta \rightarrow \infty} \pr\big( R(S) \leq \delta \big) = 1$.  $\Halmos$
\endproof\ \\

\subsection{Proof of Lemma\ \ref{brownprops00}.(\ref{supdrift11b}).}
\proof{Proof:}[Proof of Lemma\ \ref{brownprops00}.(\ref{supdrift11b})]
Precise closed-form expressions are generally well-known for the distribution of the running maximum of B.m. (see for example \citet{boukai1990explicit}).  In particular, for all $\epsilon \in (0,1)$ and $T,\gamma > 0$,
\begin{equation}\label{runmax1}
\pr\bigg(\sup_{0 \leq t \leq T}\big( {\mathcal B}^0(t) + \gamma t \big) \leq \epsilon \gamma T \bigg) = \Phi^c\big( (1 - \epsilon) \gamma T^{\frac{1}{2}} \big) - \exp(2 \epsilon \gamma^2 T) \Phi^c\big( (1 + \epsilon) \gamma T^{\frac{1}{2}} \big).
\end{equation}
Combining with Lemma\ \ref{gaussbound00}, the fact that $2 \epsilon \gamma^2 T - \frac{1}{2} (1 + \epsilon)^2 \gamma^2 T = - \frac{1}{2} (1 - \epsilon)^2 \gamma^2 T$, and the fact that our assumptions imply $\frac{1}{(1 - \epsilon)^3 \gamma^2 T} \leq \frac{\epsilon}{1 - \epsilon^2}$, we conclude that the left-hand-side of (\ref{runmax1}) is at least 
\begin{eqnarray*}
\ &\ &\ (2 \pi)^{-\frac{1}{2}} \exp\big(- \frac{1}{2} (1 - \epsilon)^2 \gamma^2 T \big) \bigg( \frac{ (1 - \epsilon) \gamma T^{\frac{1}{2}}}{\big( (1 - \epsilon) \gamma T^{\frac{1}{2}} \big)^2 + 1} - \frac{1}{(1 + \epsilon) \gamma T^{\frac{1}{2}}} \bigg)
\\&\ &\ \ \ \ \ \ \geq\ \ \ (2 \pi)^{-\frac{1}{2}} \exp\big(- \frac{1}{2} (1 - \epsilon)^2 \gamma^2 T \big) \bigg( \frac{1}{(1 - \epsilon) \gamma T^{\frac{1}{2}}} -  \frac{1}{\big((1 - \epsilon) \gamma T^{\frac{1}{2}}\big)^3} - \frac{1}{(1 + \epsilon) \gamma T^{\frac{1}{2}}} \bigg)
\\&\ &\ \ \ \ \ \ =\ \ \ (2 \pi)^{-\frac{1}{2}} \frac{1}{\gamma T^{\frac{1}{2}}} \exp\big(- \frac{1}{2} (1 - \epsilon)^2 \gamma^2 T \big) \bigg( \frac{2 \epsilon}{1 - \epsilon^2} - \frac{1}{(1 - \epsilon)^3 \gamma^2 T} \bigg)
\\&\ &\ \ \ \ \ \ \geq\ \ \ (2 \pi)^{-\frac{1}{2}} \frac{\epsilon}{1 - \epsilon^2} \frac{1}{\gamma T^{\frac{1}{2}}} \exp\big(- \frac{1}{2} (1 - \epsilon)^2 \gamma^2 T \big),
\end{eqnarray*}
completing the proof.  $\Halmos$.
\endproof

\subsection{Proof of Lemma\ \ref{brownprops00}.(\ref{brownbig1}).}\label{besselproofsec}
Although there is a vast literature on B.m. and its hitting times to various boundaries, it seems there is considerably less work on how these hitting times behave for B.m. with negative drift as the magnitude of the drift goes to zero, and that no explicit result in the literature provides the precise results which we will need in our analysis.  The proof will actually be somewhat subtle, and involve a careful analysis of the hitting times of B.m. under certain conditionings.  First we will review several properties of the so-called Three-dimensional (3-D) Bessel  Process, which will be critical for describing the relevant behaviors of conditioned B.m.
\subsubsection{The three-dimensional Bessel process.}
For any $b > 0$, let $\lbrace {\mathcal S}^b_i(t), i \geq 1 \rbrace$ denote a collection of mutually independent so-called 3-D Bessel processes initialized to $b$ (see \citet{GYor.03}).  As a notational convenience, we occasionally denote ${\mathcal S}^b_1$ by ${\mathcal S}^b$.  We now formally define ${\mathcal S}^b$ as the solution to a certain stochastic integral equation.  The stochastic integral equation
\begin{equation}\label{stochint1}
X_t = b^2 + 3 t + 2 \int_0^t |X_t|^{\frac{1}{2}} dB_s
\end{equation}
has a unique strong solution ${\mathcal X}^{b^2}(t)$, which is non-negative; we refer the reader to the survey paper of \citet{GYor.03} for details.  Then
${\mathcal S}^b$, the 3-D Bessel process initialized to $b$, is defined as $\big( {\mathcal X}^{b^2} \big)^{\frac{1}{2}}$.  The 3-D Bessel processs will be useful in our analysis, since it has the same distribution as a B.m. conditioned to hit one level before another, an object which will arise when bounding the probability of certain events.  In particular, the following is  
proven by \citet{Williams.74}, and restated by \citet{Pitman.75} Proposition 1.1.
\begin{lemma}[\citet{Pitman.75}]\label{pitman}
For any fixed $0 < b < c < \infty$, the conditional distribution of the r.v.
$$\tau^c_{ {\mathcal B}^b }\ \ \textrm{given}\ \ \big\lbrace \tau^c_{ {\mathcal B}^b } < \tau^0_{ {\mathcal B}^b } \big\rbrace$$
is identical to the distribution of the r.v. $\tau^c_{ {\mathcal S}^b }.$
Also, the conditional distribution of the process 
$${\mathcal B}^b(t)_{0 \leq t \leq \tau^c_{ {\mathcal B}^b }}\ \ \textrm{given}\ \ \big\lbrace \tau^c_{ {\mathcal B}^b } < \tau^0_{ {\mathcal B}^b } \big\rbrace$$
is identical to the distribution of the process ${\mathcal S}^b(t)_{0 \leq t \leq \tau^c_{ {\mathcal S}^b }}.$
\end{lemma}
We will also need several other generally well-known technical results w.r.t. ${\mathcal S}^b$, for use in our analysis.  
\begin{lemma}[\citet{Williams.74}]\label{bessel03bm}
${\mathcal S}^0$, namely the 3-D Bessel process initialized to 0, is distributed (on the process-level) as $\big(\sum_{i=1}^3 ({\mathcal B}^0_i)^2\big)^{\frac{1}{2}}$, where $\lbrace {\mathcal B}^0_i, i =1,2,3 \rbrace$ are independent B.m. initialized to 0 -
i.e. the radial distance process of a 3-D B.m.
\end{lemma}
More generally, an elegant construction for ${\mathcal S}^b$ (for general $b > 0$) is given by \citet{Williams.74}, where it is shown that for $b > 0$, ${\mathcal S}^b$ is distributed as the `gluing together' of two B.m.s initialized to $b$, and a 3-D Bessel process initialized to 0.  We now make this more precise.  Let $\lbrace U_b, b \geq 0 \rbrace$ denote a set of independent uniformly distributed r.v.s, where $U_b$ has the uniform distribution on $[0,b]$.  Suppose ${\mathcal B}^b_1, {\mathcal B}^b_2, {\mathcal S}^0_1, U_b$ are mutually independent and constructed on a common probability space.  Then the following is proven in \citet{Williams.74} Theorem 3.1.
\begin{lemma}[\citet{Williams.74}]\label{williams1}
For $b > 0$, define 
\begin{align*}
{\mathcal X}(t) &\stackrel{\Delta}{=} \begin{cases} {\mathcal B}^b_1(t) &\ 0 \leq t < \tau^{U_b}_{{\mathcal B}^b_1};
\\ {\mathcal B}^b_2\big( \tau^{U_b}_{{\mathcal B}^b_1} + \tau^{U_b}_{{\mathcal B}^b_2} - t \big) &\ \tau^{U_b}_{{\mathcal B}^b_1} \leq t <  \tau^{U_b}_{{\mathcal B}^b_1} + \tau^{U_b}_{{\mathcal B}^b_2};
\\ {\mathcal S}^0_1\big( t - \tau^{U_b}_{{\mathcal B}^b_1} - \tau^{U_b}_{{\mathcal B}^b_2} \big) + b &\ \tau^{U_b}_{{\mathcal B}^b_1} + \tau^{U_b}_{{\mathcal B}^b_2} \leq t < \infty.
\end{cases}
\end{align*}
Then the distribution of the process ${\mathcal X}$ is identical to the distribution of the process ${\mathcal S}^b_1$.
\end{lemma}
\subsubsection{Preliminary results on hitting times of the three-dimensional Bessel process.}
We will also need several additional preliminary results regarding the hitting times of the 3-D Bessel process.  Although these hitting times are generally well-studied (see \citet{Kent.78}, \citet{GS.79}, \citet{PY.81}, \citet{BR.06}, \citet{BMR.10}), we include proofs with associated simple explicit bounds, both for completeness, and as the associated results in the literature do not seem to be in a form amenable to our particular application.  The first result we will need proves that the time it takes ${\mathcal S}^b$ to hit a high level is sufficiently large with sufficiently high probability.  We defer the proof to the end of Subsection\ \ref{besselproofsec}.
\begin{lemma}\label{besselhit1}
For all $M > b > 0$ and $x > 0$, 
$$\pr\big(\tau^M_{{\mathcal S}^b} > x \big) \geq 1 - \frac{2 b}{M} - \frac{100 x}{(M - b)^2}.$$
\end{lemma}

We will also need a more subtle law-of-the-iterated-logarithm type result regarding the probability that a 3-D Bessel process never goes below a certain boundary.  Although related results are generally well-known (see \citet{Hkk.03}), we now state an explicit variant customized to the needs of our own proof.  We again defer the proof to the end of Subsection\ \ref{besselproofsec}.
\begin{lemma}\label{itlog2}
For all $C \geq 1$ and $b \geq 4 C^2$,
$$\pr\bigg( \inf_{t \geq 0}  \big( {\mathcal S}^b(t) - C \log^{\frac{1}{2}}(3 + 3 t) \big) > 0\bigg) \geq 1 - 10^6 \times C^5 \times b^{-\frac{1}{2}}.$$
Namely, the probability that ${\mathcal S}^b(t)$ stays above $C \log^{\frac{1}{2}}(3 + 3 t)$ for all times $t$ is as least $1 - 10^6 \times C^5 \times b^{-\frac{1}{2}}$.
\end{lemma}

\subsubsection{Proof of  Lemma\ \ref{brownprops00}.(\ref{brownbig1}).}
We now complete the proof of Lemma\ \ref{brownprops00}.(\ref{brownbig1}).
\proof{Proof:}[Proof of Lemma\ \ref{brownprops00}.(\ref{brownbig1})]
Let $\zeta \stackrel{\Delta}{=} 10^3 \times C$, 
$${\mathcal E}_0 \stackrel{\Delta}{=} \bigg\lbrace \sup_{t \geq 0} \bigg( {\mathcal B}^{-b}(t) + C \log^{\frac{1}{2}}(3 + 3 t) - \gamma t\bigg)\ \leq 0 \bigg\rbrace\ \ \ ,\ \ \ 
{\mathcal E}_1 \stackrel{\Delta}{=}  \bigg\lbrace \tau^{-\frac{\zeta}{\gamma}}_{{\mathcal B}^{-b}} < \tau^{0}_{{\mathcal B}^{-b}} \bigg\rbrace.$$
In that case, it follows from Lemma\ \ref{brownprops00}.(\ref{hitwhich}) and the symmetries of B.m. that 
\begin{equation}\label{bzero1}
\pr({\mathcal E}_0) = \frac{b \gamma}{\zeta} \pr({\mathcal E}_0 | {\mathcal E}_1).
\end{equation}
Suppose the B.m. ${\mathcal B}^0$ and the 3-D Bessel process ${\mathcal S}^{b}$ are independent, and constructed on a common probability space.  Let us define
$${\mathcal E}_2 \stackrel{\Delta}{=} \bigg\lbrace \sup_{0 \leq t \leq \tau^{\frac{\zeta}{\gamma}}_{{\mathcal S}^{b}}} \bigg( - {\mathcal S}^{b}(t) + C \log^{\frac{1}{2}}(3 + 3 t) - \gamma t\bigg)\ \leq 0 \bigg\rbrace,$$
$${\mathcal E}_3 \stackrel{\Delta}{=} \bigg\lbrace \sup_{t \geq 0} \bigg( - \frac{\zeta}{\gamma} + {\mathcal B}^0(t) + C \log^{\frac{1}{2}}(3 + 3 \tau^{\frac{\zeta}{\gamma}}_{{\mathcal S}^{b}} + 3 t) - \gamma (t + \tau^{\frac{\zeta}{\gamma}}_{{\mathcal S}^{b}})\bigg)\ \leq 0 \bigg\rbrace.$$
It then follows from Lemma\ \ref{pitman}, the symmetries of B.m., and the independent increments and strong Markov properties of B.m., that 
\begin{equation}\label{bzero2aaa}
\pr({\mathcal E}_0 | {\mathcal E}_1) = \pr({\mathcal E}_2 , {\mathcal E}_3).
\end{equation}
We now bound $\pr({\mathcal  E}_2)$ and $\pr({\mathcal E}_3)$, beginning with $\pr({\mathcal E}_2)$.
Let
$${\mathcal E}_4 \stackrel{\Delta}{=} \bigg\lbrace \inf_{t \geq 0} \bigg( {\mathcal S}^{b}(t) - C \log^{\frac{1}{2}}(3 + 3 t) \bigg)\ \geq 0 \bigg\rbrace,$$
and note that by the natural monotonicities of the relevant quantities and Lemma\ \ref{itlog2},
\begin{equation}\label{bzero2}
\pr({\mathcal E}_2)\ \ \ \geq\ \ \ \pr({\mathcal E}_4)\ \ \ \geq\ \ \ 1 - 10^6 \times C^5 \times b^{-\frac{1}{2}}.
\end{equation}
We next bound $\pr({\mathcal E}_3)$.  Let 
$${\mathcal E}_5 \stackrel{\Delta}{=}  \bigg\lbrace \tau^{\frac{\zeta}{\gamma}}_{{\mathcal S}^{b}} > \frac{16 C^2}{\gamma^2} \bigg\rbrace\ \ \ ,\ \ \ 
{\mathcal E}_6 \stackrel{\Delta}{=} \bigg\lbrace \sup_{t \geq 0} \bigg( {\mathcal B}^0(t) - \frac{\gamma}{2} t \bigg) \leq \frac{\zeta}{\gamma} \bigg\rbrace.$$
It follows from straightforward algebra that $C \log^{\frac{1}{2}}(3 + 3 y) < \frac{\gamma}{2} y$ for all $y > \frac{16 C^2}{\gamma^2}$.  Indeed, as $(a + b)^{\frac{1}{2}} \leq a^{\frac{1}{2}} + b^{\frac{1}{2}}$ for all $a,b > 0$, we find that 
$\log^{\frac{1}{2}}(3 + 3 y) \leq \log^{\frac{1}{2}}(3) + \log^{\frac{1}{2}}(y+1)$.  As the exponential inequality implies that $\log(y + 1) \leq y$, and our assumptions imply that $\log^{\frac{1}{2}}(3) < y^{\frac{1}{2}}$, the desired inequality follows.
We may thus conclude from the natural monotonicities of the relevant quantities that
\begin{eqnarray*}
\pr({\mathcal E}_3 | {\mathcal E}_5) &\geq& \pr\Bigg( \sup_{t \geq 0} \bigg( - \frac{\zeta}{\gamma} + {\mathcal B}^0(t) + \frac{\gamma}{2}( \tau^{\frac{\zeta}{\gamma}}_{{\mathcal S}^{b}} + t) - \gamma (t + \tau^{\frac{\zeta}{\gamma}}_{{\mathcal S}^{b}})\bigg)\ \leq 0 \Bigg) 
\\&=& \pr\Bigg( \sup_{t \geq 0} \bigg( - \frac{\zeta}{\gamma} + {\mathcal B}^0(t)  - \frac{\gamma}{2} t - \frac{\gamma}{2} \tau^{\frac{\zeta}{\gamma}}_{{\mathcal S}^{b}} \bigg)\ \leq 0 \Bigg) 
\\&\geq& \pr\bigg( \sup_{t \geq 0} \big( {\mathcal B}^0(t)  - \frac{\gamma}{2} t \big)\ \leq \frac{\zeta}{\gamma} \bigg)\ \ \ =\ \ \ \pr({\mathcal E}_6),
\end{eqnarray*}
and we conclude that
\begin{equation}\label{bzero3}
\pr({\mathcal E}_3)\ \ \ \geq\ \ \ \pr({\mathcal E}_5) \times \pr({\mathcal E}_6).
\end{equation}
Lemma\ \ref{besselhit1}, combined with our assumptions which imply that $b < \frac{\zeta}{4 \gamma}$ and some straightforward algebra, ensures that
$$
\pr({\mathcal E}_5)\ \ \geq\ \ \frac{1}{2} - \frac{6400 C^2}{\zeta^2}\ \ =\ \ \frac{1}{2} - \frac{6400}{10^6}.
$$
As it follows from Lemma\ \ref{brownprops00}.(\ref{supdrift1}) that
$$\pr({\mathcal E}_6)\ \ \ =\ \ \ 1 - \exp(- \zeta)\ \ \ \geq\ \ \ 1 - \exp(- 10^3),
$$
we conclude from (\ref{bzero3}) and some straightforward algebra that 
\begin{equation}\label{bzer00a}
\pr({\mathcal E}_3)\ \ \ \geq\ \ \ (\frac{1}{2} - \frac{6400}{10^6}) \times \big(1 - \exp(- 10^3) \big)\ \ \ \geq\ \ \ .49.
\end{equation}
Combining (\ref{bzero2}) and (\ref{bzer00a}) with (\ref{bzero2aaa}) and a union bound, we conclude that $\pr({\mathcal E}_0 | {\mathcal E}_1) \geq .49 - 10^6 \times C^5 \times b^{-\frac{1}{2}}$.
As our assumptions imply that $b \geq 10^{13} \times C^{10}$, it follows from some straightforward algebra that $\pr({\mathcal E}_0 | {\mathcal E}_1) \geq .17$.
Combining with (\ref{bzero1}) completes the proof.  $\Halmos$.
\endproof

\subsubsection{Proof of Lemma\ \ref{besselhit1}.}
\proof{Proof:}[Proof of Lemma\ \ref{besselhit1}]
It follows from Lemma\ \ref{williams1} and a union bound that $\pr\big(\tau^M_{{\mathcal S}^b} > x \big)$ is at least 
$$\pr\big(\tau^{U_b}_{{\mathcal B}^b_1} < \tau^{M}_{{\mathcal B}^b_1}\ \ ,\ \ \tau^{U_b}_{{\mathcal B}^b_2} < \tau^{M}_{{\mathcal B}^b_2}\big) \times \pr\big(\tau^{M-b}_{{\mathcal S}^0} > x\big),
$$
which is itself at least
$$\pr\big(\tau^{0}_{{\mathcal B}^b_1} < \tau^{M}_{{\mathcal B}^b_1}\big)^2 \times \pr\big(\tau^{M-b}_{{\mathcal S}^0} > x\big).
$$
Further applying Lemma\ \ref{brownprops00}.(\ref{hitwhich}) and the relevant symmetries of B.m., we find that 
\begin{equation}\label{bessel002}
\pr\big(\tau^M_{{\mathcal S}^b} > x \big) \geq (1 - \frac{b}{M})^2 \times \pr\big(\tau^{M-b}_{{\mathcal S}^0} > x\big).
\end{equation}
It follows from Lemma\ \ref{bessel03bm}, Lemma\ \ref{brownprops00}.(\ref{supform}), a union bound, and the symmetries of B.m. that
\begin{eqnarray}
\pr\big(\tau^{M-b}_{{\mathcal S}^0} > x\big) &=& \pr\big(\tau^{(M-b)^2}_{\sum_{i=1}^3 ({\mathcal B}^0_i)^2} > x\big) \nonumber
\\&\geq& \prod_{i=1}^3 \pr\bigg(\tau^{\frac{(M-b)^2}{3}}_{({\mathcal B}^0_i)^2} > x\bigg) \nonumber
\\&=& \pr^3\bigg(\tau^{3^{-\frac{1}{2}}(M-b)}_{|{\mathcal B}^0_1|} > x\bigg) \nonumber
\\&\geq& \bigg(1 - 2 \pr\big(\tau^{3^{-\frac{1}{2}}(M-b)}_{{\mathcal B}^0_1} \leq x\big) \bigg)^3 \nonumber
\\&=& \bigg(1 - 2 \pr\big( \sup_{0 \leq s \leq x} {\mathcal B}^0_1(s) > 3^{-\frac{1}{2}}(M-b) \big) \bigg)^3 \nonumber
\\&=& \bigg(1 - 4 \Phi^c\big( (3 x)^{-\frac{1}{2}}(M-b) \big) \bigg)^3.\label{bbbmmm1}
\end{eqnarray}
Combining the above with Lemma\ \ref{gaussbound00}, the fact that $\exp(x) > x$ for all $x > 0$, the fact that $(1 - a)^k \geq 1 - k a$ for all $a \in (0,1)$ and $k \geq 1$, the fact that $(1 - a)(1-b) \geq 1 - a - b$ for all $a,b \in (0,1)$, 
and some straightforward algebra completes the proof.  $\Halmos$.
\endproof

\subsubsection{Proof of Lemma\ \ref{itlog2}.}
\proof{Proof:}[Proof of Lemma\ \ref{itlog2}]
Suppose that ${\mathcal B}^b_1, {\mathcal B}^b_2, U_b$, and ${\mathcal S}^0_1$ have been used to construct ${\mathcal S}^b_1$ as dictated by the construction given in Lemma\ \ref{williams1}.
Let us define the events 
$${\mathcal E}_0 \stackrel{\Delta}{=}  \big\lbrace \inf_{t \geq 0}  \big( {\mathcal S}^b_1(t) - C \log^{\frac{1}{2}}(3 + 3 t) \big) > 0 \big\rbrace,$$
$${\mathcal E}_1 \stackrel{\Delta}{=} \big\lbrace U_b > C \log^{\frac{1}{2}}(3 + 3 \tau^{U_b}_{ {\mathcal B}^b_1 } + 3 \tau^{U_b}_{ {\mathcal B}^b_2 }) \big\rbrace,$$
$${\mathcal E}_2 \stackrel{\Delta}{=} \Bigg\lbrace \inf_{t \geq 0} \bigg( b + {\mathcal S}^0_1(t) -  C \log^{\frac{1}{2}}(3 + 3 \tau^{U_b}_{ {\mathcal B}^b_1 } + 3 \tau^{U_b}_{ {\mathcal B}^b_2 } + 3 t) \bigg) > 0 \Bigg\rbrace.$$
In that case, it may be easily verified using Lemma\ \ref{williams1} and a union bound that 
\begin{equation}\label{bessliluser1}
\pr({\mathcal E}_0)\ \ \ \geq\ \ \ \pr({\mathcal E}_1, {\mathcal E}_2)\ \ \ \geq\ \ \ 1 - \pr({\mathcal E}^c_1) - \pr({\mathcal E}^c_2).
\end{equation}
We proceed by bounding $\pr({\mathcal E}^c_1)$ and $\pr({\mathcal E}^c_2)$, and begin by bounding $\pr({\mathcal E}^c_1)$.  By conditioning on $U_b$, we find that
\begin{equation}\label{bessliluser2}
\pr({\mathcal E}^c_1) = b^{-1} \int_0^b \pr\bigg( \tau^{u}_{ {\mathcal B}^b_1 } + \tau^{u}_{ {\mathcal B}^b_2 } \geq \frac{1}{3} \exp\big( (\frac{u}{C})^2\big) - 1 \bigg) du.
\end{equation}
As $\tau^{u}_{ {\mathcal B}^b_1 }$ and $\tau^{u}_{ {\mathcal B}^b_2 }$ are i.i.d., it follows follows from (\ref{bessliluser2}), Lemma\ \ref{brownprops00}.(\ref{hitdist00}), a union bound, and the fact that our assumptions imply that  
$$\frac{1}{3} \exp\big( (\frac{u}{C})^2\big) - 1 > \frac{1}{6} \exp\big( (\frac{u}{C})^2\big)\ \ \textrm{and}\ \ \frac{u}{C} > 1\ \ \textrm{for all}\ \ u > b^{\frac{1}{2}},
$$
that
\begin{eqnarray}
\pr({\mathcal E}^c_1) &\leq& b^{-1} \int_0^{b^{\frac{1}{2}}} 1 du + b^{-1} \int_{b^{\frac{1}{2}}}^b \pr\bigg( \tau^{u}_{ {\mathcal B}^b_1 } + \tau^{u}_{ {\mathcal B}^b_2 } \geq \frac{1}{6} \exp\big( (\frac{u}{C})^2\big) \bigg) du \nonumber
\\&\leq& b^{-1} \int_0^{b^{\frac{1}{2}}} 1 du + 2 b^{-1} \int_{b^{\frac{1}{2}}}^b \pr\bigg( \tau^{u}_{ {\mathcal B}^b_1 } \geq \frac{1}{12} \exp\big( (\frac{u}{C})^2\big) \bigg) du \nonumber
\\&\leq& b^{-\frac{1}{2}} + 2 b^{-1} \int_{b^{\frac{1}{2}}}^b (b - u) \bigg( \frac{1}{12} \exp\big( (\frac{u}{C})^2\big) \bigg)^{-\frac{1}{2}} du \nonumber
\\&\leq& b^{-\frac{1}{2}} + 8  \int_{b^{\frac{1}{2}}}^{\infty} \frac{u}{C} \exp\big( - \frac{1}{2} (\frac{u}{C})^2\big) du\ \ \ =\ \ \ b^{-\frac{1}{2}} + 8 C \exp( - \frac{b}{2 C^2} ). \label{bessliluser3}
\end{eqnarray}
We next bound $\pr({\mathcal E}^c_2)$, which by a union bound is at most
\begin{eqnarray}
\ &\ &\ \ \ \pr\bigg(\tau^{U_b}_{ {\mathcal B}^b_1 } + \tau^{U_b}_{ {\mathcal B}^b_2 } \geq \frac{1}{3} \exp(\frac{b^2}{4 C}) - \frac{4}{3}\bigg)  \label{bessliluser4}
\\&\ &\ \ \ \ \ \ \ +\ \ \ \pr\Bigg( \inf_{t \geq 0} \bigg( b + {\mathcal S}^0_1(t) -  C \log^{\frac{1}{2}}\big( \exp( \frac{b^2}{4 C}) - 1+ 3 t\big) \leq 0 \bigg) \Bigg). \label{bessliluser5}
\end{eqnarray}
Applying an argument nearly identical to that used to bound $\pr({\mathcal E}^c_1)$, combined with the fact that our assumptions imply that $\frac{1}{3}\exp(\frac{b^2}{4 C}) - \frac{4}{3} \geq \frac{1}{6} \exp( \frac{b^2}{4 C} )$, we find that (\ref{bessliluser4}) is at most 
\begin{eqnarray}
\ &\ &\ \ \ b^{-1} \int_0^b 2 \pr\big( \tau^{u}_{ {\mathcal B}^b_1 } \geq \frac{1}{12} \exp( \frac{b^2}{4 C^2}) \big) du \nonumber
\\&\ &\ \ \ \ \ \ \leq\ \ \ 8 b^{-1} \int_0^b (b - u) \exp( - \frac{b^2}{8 C^2}) du \ \ \ \leq\ \ \ 8 b \exp( - \frac{b^2}{8 C^2} ).\label{bessliluser6}
\end{eqnarray}
As it is easily verified that $\log^{\frac{1}{2}}(x + y) \leq \log^{\frac{1}{2}}(x + 1) + \log^{\frac{1}{2}}(y + 1)$ for all $x,y > 0$, it follows that 
$$
C \log^{\frac{1}{2}}\big( \exp( \frac{b^2}{4 C^2} ) - 1+ 3 t \big) \leq \frac{b}{2} + C \log^{\frac{1}{2}}(1 + 3 t),$$
and thus (\ref{bessliluser5}) is at most 
\begin{equation}\label{bessliluser7}
\pr\Bigg( \inf_{t \geq 0} \bigg( \frac{b}{2}+ {\mathcal S}^0_1(t) -  C \log^{\frac{1}{2}}(1 + 3 t) \leq 0 \bigg) \Bigg).
\end{equation}
Let $T' \stackrel{\Delta}{=} \frac{1}{3}\big(\exp( \frac{b^2}{4 C^2} ) - 1 \big)$ denote the smallest time $t$ for which $C \log^{\frac{1}{2}}(1 + 3 t) \geq \frac{b}{2}$, where we note that our assumptions imply $\lfloor T' \rfloor \geq \frac{1}{6} \exp( \frac{b^2}{4 C^2} )$.
In that case, it follows from a straightforward contradiction argument and a union bound that 
(\ref{bessliluser7}) is at most 
\begin{equation}\label{bessliluser8}
\sum_{k=\lfloor T' \rfloor}^{\infty} \pr\bigg( \inf_{k \leq t \leq k + 1} {\mathcal S}^0_1(t)  \leq C \log^{\frac{1}{2}}(4 + 3 k) \bigg).
\end{equation}
By Lemma\ \ref{bessel03bm}, a union bound, the independence of $\lbrace {\mathcal B}^0_i, i = 1,2,3 \rbrace$, and the stationary increments property of B.m., we find that for all $k \geq 0$, $\pr\bigg( \inf_{k \leq t \leq k + 1} {\mathcal S}^0_1(t)  \leq C \log^{\frac{1}{2}}(4 + 3 k) \bigg)$ is at most 
\begin{eqnarray}
\ &\ &\ \ \pr\bigg( \inf_{k \leq t \leq k + 1} \sum_{i=1}^3 \big({\mathcal B}^0_i(t)\big)^2  \leq C^2 \log(4 + 3 k) \bigg) \nonumber
\\&\ &\ \ \ \ \leq\ \ \ \pr\bigg( \bigcap_{i=1}^3 \big\lbrace \inf_{k \leq t \leq k + 1} \big|{\mathcal B}^0_i(t)|  \leq C \log^{\frac{1}{2}}(4 + 3 k) \big\rbrace \bigg) \nonumber
\\&\ &\ \ \ \ =\ \ \ \pr^3\bigg( \inf_{k \leq t \leq k + 1} |{\mathcal B}^0_1(t)|  \leq C \log^{\frac{1}{2}}(4 + 3 k) \bigg) \nonumber
\\&\ &\ \ \ \ \leq\ \ \ \Bigg(\pr\big(|{\mathcal B}^0_1(k)| \leq 3 C \log^{\frac{1}{2}}(4 + 3 k)\big) + \pr\big( \sup_{0 \leq t \leq 1} |{\mathcal B}^0_1(t)| \geq 2 C \log^{\frac{1}{2}}(4 + 3 k)\big) \Bigg)^3. \label{bessliluser9}
\end{eqnarray}
We next bound $\pr\big(|{\mathcal B}^0_1(k)| \leq 3 C \log^{\frac{1}{2}}(4 + 3 k)\big)$.  Let $N$ denote a normally distributed r.v. with mean 0 and unit variance.  Applying the basic properties of B.m. with the fact that $\phi(x) \leq (2 \pi)^{-\frac{1}{2}}$ for all $x \in \bbr$ (i.e. the Gaussian density is bounded by $(2 \pi)^{-\frac{1}{2}}$), as well as the easily verified fact that $\log^{\frac{1}{2}}(4 + 3 k) \leq 2 k^{\frac{1}{8}}$ for all $k \geq 1$, we conclude that
\begin{equation}\label{bessliluser10}
\pr\big(|{\mathcal B}^0_1(k)| \leq 3 C \log^{\frac{1}{2}}(4 + 3 k)\big)\ \ \ =\ \ \ \pr\big(|N| \leq  3 C k^{-\frac{1}{2}} \log^{\frac{1}{2}}(4 + 3 k)\big)\ \ \ \leq\ \ \ 6 C k^{-\frac{3}{8}}. 
\end{equation}
Furthermore, by Lemma\ \ref{brownprops00}.(\ref{supform}), Lemma\ \ref{gaussbound00}, a union bound, and the fact that by assumption $C \geq 1$,
\begin{equation}\label{bessliluser11}
\pr\big( \sup_{0 \leq t \leq 1} |{\mathcal B}^0_1(t)| \geq 2 C \log^{\frac{1}{2}}(4 + 3 k)\big)\ \ \ \leq\ \ \ \frac{4}{(2 \pi)^{\frac{1}{2}}} \exp\bigg( - \frac{1}{2}\big(2 C \log^{\frac{1}{2}}(4 + 3 k)\big)^2\bigg) \leq 2 C k^{-\frac{3}{8}}.
\end{equation}
Combining (\ref{bessliluser7}) - (\ref{bessliluser11}), it follows that (\ref{bessliluser5}) is at most
\begin{eqnarray}
\sum_{k = \lfloor T' \rfloor}^{\infty} (8 C k^{-\frac{3}{8}})^3 &\leq& 600 C^3 \sum_{k = \lfloor T' \rfloor}^{\infty} k^{-\frac{9}{8}}\nonumber
\\&\leq& 600 C^3 \int_{ \frac{1}{6} \exp( \frac{b^2}{4 C^2} ) }^{\infty} u^{-\frac{9}{8}} du\ \ \ \leq\ \ \ 10^4 C^3 \exp( - \frac{b^2}{32 C^2} ). \label{bessliluser12}
\end{eqnarray}
Combining (\ref{bessliluser1}), (\ref{bessliluser3}), (\ref{bessliluser6}), and (\ref{bessliluser12}), we conclude that 
\begin{equation}\label{bessliluser13}
\pr({\mathcal E}_0) \geq 1 - b^{-\frac{1}{2}} - 8 C \exp( - \frac{b}{2 C^2} ) - 8 b \exp( - \frac{b^2}{8 C^2} ) - 10^4 C^3 \exp( - \frac{b^2}{32 C^2} ).
\end{equation}
Applying the fact that $\exp(x) > x$ for all $x > 0$ and some straightforward algebra completes the proof.  $\Halmos$.
\endproof

\subsection{Proof of Lemma\ \ref{ousup2}.}
\proof{Proof:}[Proof of Lemma\ \ref{ousup2}]
By stationarity and a union bound, it suffices to demonstrate that
$$\sum_{k=0}^{\infty} \pr\big( \sup_{0 \leq t \leq 1} |{\mathcal U}^{\rho}(t)| \geq 6 \exp(\rho) \log^{\frac{1}{2}}(3 + 3 k) \big) \leq \frac{1}{2} .
$$
It is well-known that one can construct ${\mathcal U}^{\rho}$ on the same probability space as a mutually independent $N(0,1)$ r.v. $N$ and B.m. ${\mathcal B}^0$, s.t. 
${\mathcal U}^{\rho}(t) = \exp(- \rho t) \bigg( N + {\mathcal B}^0\big(\exp(2 \rho t) - 1\big) \bigg)$ for all $t \geq 0$.  Hence on the same probability space, 
$$\pr\bigg( \sup_{t \geq 0} \frac{|{\mathcal U}^{\rho}(t)|}{|N| + |{\mathcal B}^0\big(\exp(2 \rho t) - 1\big)|} \leq 1 \bigg) = 1.$$
Namely, with probability one, $|{\mathcal U}^{\rho}(t)|$ is (for all times $t \geq 0$) dominated by $|N| + |{\mathcal B}^0\big(\exp(2 \rho t) - 1\big)|$.
To complete the proof it thus suffices to demonstrate that 
\begin{equation}\label{showoou2}
\sum_{k=0}^{\infty} \pr\big( |N| + \sup_{0 \leq t \leq \exp(2 \rho) - 1} |{\mathcal B}^0(t)| \geq 6 \exp(\rho) \log^{\frac{1}{2}}(3 + 3 k) \big) \leq \frac{1}{2}.
\end{equation}
Applying a union bound, we conclude that the left-hand-side of (\ref{showoou2}) is at most
\begin{eqnarray*}
&\ &\ \sum_{k=0}^{\infty} \pr\big( |N| \geq 2 \exp(\rho) \log^{\frac{1}{2}}(3 + 3 k) \big) + \sum_{k=0}^{\infty} \pr\big( \sup_{0 \leq t \leq \exp(2 \rho) - 1} |{\mathcal B}^0(t)| \geq 4 \exp(\rho) \log^{\frac{1}{2}}(3 + 3 k) \big)
\\&\ &\ \ \ \leq\ \ \ 2 \sum_{k=0}^{\infty} \pr\big( N \geq 2 \exp(\rho) \log^{\frac{1}{2}}(3 + 3 k) \big) + 2 \sum_{k=0}^{\infty} \pr\big( \sup_{0 \leq t \leq \exp(2 \rho) - 1} {\mathcal B}^0(t) \geq 2 \exp(\rho) \log^{\frac{1}{2}}(3 + 3 k) \big),
\end{eqnarray*}
which by applying Lemma\ \ref{brownprops00}.(\ref{supform}) and bounding $\exp(2 \rho) - 1$ by $\exp(2 \rho)$ is at most 
\begin{equation}\label{showoou3}
2 \sum_{k=0}^{\infty} \Phi^c\big(2 \exp(\rho) \log^{\frac{1}{2}}(3 + 3 k) \big) + 4 \sum_{k=0}^{\infty} \Phi^c\bigg(2 \log^{\frac{1}{2}}(3 + 3 k) \bigg).
\end{equation}
Further noting that $\exp(\rho) \geq 1$ and $\log^{\frac{1}{2}}(3 + 3 k) \geq 1$ for all $k \geq 0$, and applying Lemma\ \ref{gaussbound00}, it follows that (\ref{showoou3}) is at most
\begin{eqnarray*}
6 \sum_{k=0}^{\infty} \Phi^c\big(2 \log^{\frac{1}{2}}(3 + 3 k) \big) &\leq& \frac{6}{(2 \pi)^{\frac{1}{2}}} \sum_{k=0}^{\infty} (3 + 3 k)^{-2}
\\&=& \frac{\pi^{\frac{3}{2}}}{9 \times 2^{\frac{1}{2}}}\ \ \ <\ \ \ \frac{1}{2}.
\end{eqnarray*}
Combining the above completes the proof.  $\Halmos$.
\endproof

\subsection{Proof of Theorem\ \ref{compareprocess}, Lemma\ \ref{boundvarzztop}, Lemma\ \ref{poscovar0}, and Corollary\ \ref{boundDsupmean}.}
In this section we complete the proof of Theorem\ \ref{compareprocess}, as well as  Lemma\ \ref{boundvarzztop}, Lemma\ \ref{poscovar0}, and Corollary\ \ref{boundDsupmean}.  All these results will follow from a careful analysis of the covariance of ${\mathcal D}$, i.e. the covariance structure of an equilibrium renewal process.
\subsubsection{Background on stationary renewal processes and their covariances.}
Recall that ${\mathcal N}_1$ is an equilibrium renewal process with renewal distribution $S$, and that $N_1(t)$ denotes the number of renewals up to time $t$.  Similarly, let ${\mathcal N}_o$ denote an ordinary renewal process with renewal distribution $S$, and $N_o(t)$ the corresponding number of renewals up to time $t$.  Much is known about the covariance structure of ${\mathcal N}_1$ and ${\mathcal N}_o$, and we now review several such properties for use in our analysis.  Let $f(t) \stackrel{\Delta}{=} V[N_1(t)] - \mu c^2_S t$.  Then the following lemma follows immediately from \citet{Daley.78}, \citet{Daley.80}, and \citet{Lorden.70}.
\begin{lemma}[\citet{Daley.78,Daley.80,Lorden.70}]\label{reuserenew}
$\E[N_o(t)]$ is a monotone increasing and integrable function, and for all $t \geq 0$,
\begin{equation}\label{frep1}
f(t) = 2 \mu \int_0^t \bigg( \big( \E[N_o(s)] + 1 - \mu s \big) - \frac{1}{2}\big( 1 + c^2_S \big) \bigg) ds.
\end{equation}
Furthermore, 
\begin{equation}\label{supf0}
\sup_{t \geq 0} |f(t)| \leq \frac{4}{3} \mu^3 \E[S^3] + \frac{1}{4} \mu^4 \big(\E[S^2]\big)^2 \leq \alpha_S;
\end{equation}
and for all $s \geq 0$, 
\begin{equation}\label{inbound1}
0 \leq  \E[N_o(s)] + 1 - \mu s \leq \mu^2 \E[S^2].
\end{equation}
\end{lemma}
We now derive several useful implications of Lemma\ \ref{reuserenew}.
\begin{corollary}\label{lipandrest}
$f$ is a Lipschitz-continuous function with Lipschitz constant less than $\mu^3 \E[S^2] \leq \alpha_S$.  Furthermore, for all $0 \leq s \leq t$:
\begin{enumerate}[(i)]
\item $\E[{\mathcal D}(s) {\mathcal D}(t)] = \mu c^2_S s + \frac{1}{2}\big( f(s) + f(t) - f(t-s)\big)$; \label{lipb1}
\item $\E[{\mathcal D}(s) {\mathcal D}(t)]  \geq 0$; \label{lipb2}
\item $\big| \E[{\mathcal D}(s) {\mathcal D}(t)] - \E[{\mathcal D}^2(s)] \big| \leq \mu^3 \E[S^2] (t - s) \leq \alpha_S (t - s)$; \label{lipb3}
\item $\big| \E[{\mathcal D}(s) {\mathcal D}(t)] - \E[{\mathcal D}^2(s)] \big| \leq 2 \mu^3 \E[S^3] + \frac{3}{8} \mu^4 \big(\E[S^2]\big)^2 \leq \alpha_S$. \label{lipb4}
\end{enumerate}
\end{corollary}
\proof{Proof:}
That $f$ is Lipschitz with Lipshitz constant at most $\mu^3 \E[S^2]$ follows from (\ref{frep1}), (\ref{inbound1}), and straightforward algebra.
\\\\We next prove (\ref{lipb1}).  Note that 
$$\E[{\mathcal D}(s) {\mathcal D}(t)] = -\frac{1}{2} \E[ \big({\mathcal D}(t) - \mathcal{D}(s) \big)^2 ] + \frac{1}{2} \E[{\mathcal D}^2(t)] + \frac{1}{2} \E[{\mathcal D}^2(s)].$$
Combining with the stationary increments property and definitions completes the proof.
\\\\We now prove (\ref{lipb2}).  Let $\lbrace X_k, k \geq 1 \rbrace$ denote the ordered sequence of renewal intervals in process ${\mathcal N}_1$.  Then from definitions,
\begin{equation}\label{poscor1}
V[N_1(s) , N_1(t)] = \sum_{i=1}^{\infty} \sum_{j=1}^{\infty} \E[ I( \sum_{k=1}^i X_k \leq s) I( \sum_{k=1}^j X_k \leq t) ]  - \mu^2 s t.
\end{equation}
Note that for $j \leq i$, one has that 
\begin{eqnarray}
\E[ I( \sum_{k=1}^i X_k \leq s) I( \sum_{k=1}^j X_k \leq t) ] &=& \E[ I( \sum_{k=1}^i X_k \leq s) ] \nonumber
\\&\geq& \E[ I( \sum_{k=1}^i X_k \leq s) ]  \E[ I( \sum_{k=1}^j X_k \leq t) ].\label{poscor2}
\end{eqnarray}
Alternatively, suppose $j \geq i+1$.  Let $Y^1 \stackrel{\Delta}{=} \sum_{k=1}^i X_k$, $Y^2 \stackrel{\Delta}{=} \sum_{k=i+1}^j X_k$, and $Y^3 \stackrel{\Delta}{=} t - Y^1$.  Then $Y^1$ and $Y^2$ are independent, $Y^2$ and $Y^3$ are independent, and \begin{eqnarray}
\E[ I( \sum_{k=1}^i X_k \leq s) I( \sum_{k=1}^j X_k \leq t) & = & \E[ I( Y^1 \leq s) I( Y^1 + Y^2 \leq t) ] \nonumber
\\&=& \E[ I( Y^3 \geq t-s ) I( Y^3 \geq Y^2) ] \nonumber
\\&=& \E[ I(Y^2 \geq t-s) I(Y^3 \geq Y^2) + I(Y^2 < t-s) I(Y^3 \geq t-s) ].\label{poscor3}
\end{eqnarray}
Let $Y^3_a,Y^3_b$ denote two r.v.s, each distributed as $Y^3$, where $Y^3_a,Y^3_b,Y^2$ are mutually independent.  Then by linearity of expectation, (\ref{poscor3}) equals
\begin{eqnarray}
&\ &\E[ I(Y^2 \geq t-s) I(Y^3_a \geq Y^2) + I(Y^2 < t-s) I(Y^3_b \geq t-s) ]\nonumber
\\&\indent&\ \ \ \ \ \ \ \geq\ \ \ \E[ I(Y^2 \geq t-s) I(Y^3_a \geq Y^2)I(Y^3_b \geq t-s)] \nonumber
\\&\indent&\ \ \ \ \ \ \ \ \ \ \ \ \ \ \ \ + \E[I(Y^2 < t-s) I(Y^3_a \geq Y^2) I(Y^3_b \geq t-s) ] \nonumber
\\&\indent&\ \ \ \ \ \ \ =\ \ \ \E[ I(Y^3_a \geq Y^2)I(Y^3_b \geq t-s) ] \nonumber
\\&\indent&\ \ \ \ \ \ \ =\ \ \ \E[ I(Y^1 + Y^2 \leq t) ]\E[ I(Y^1 \leq s) ]\nonumber.
\\&\indent&\ \ \ \ \ \ \ =\ \ \ \E[ I( \sum_{k=1}^i X_k \leq s)] \E[I( \sum_{k=1}^j X_k \leq t) ] \label{poscor4}.
\end{eqnarray}
Combining (\ref{poscor1}) - (\ref{poscor4}), we find that
\begin{eqnarray*}
V[N^e(s) , N^e(t)] &\geq& \sum_{i=1}^{\infty} \sum_{j=1}^{\infty} \E[ I( \sum_{k=1}^i X_k \leq s)] \E[I( \sum_{k=1}^j X_k \leq t) ] - \mu^2 s t
\\&=& \E\big[ \sum_{i=1}^{\infty} I( \sum_{k=1}^i X_k \leq s)\big] \E\big[\sum_{j=1}^{\infty} I( \sum_{k=1}^j X_k \leq t)\big] - \mu^2 s t
\\&=& \E[N^e(s)] \E[N^e(t)] - \mu^2 st\ \ \ =\ \ \ 0,
\end{eqnarray*}
completing the proof.
\\\\Finally, we prove both (\ref{lipb3}) and (\ref{lipb4}).  It follows from (\ref{lipb1}) that $\E[{\mathcal D}(s) {\mathcal D}(t)] - \E[{\mathcal D}^2(s)] = \frac{1}{2}\big( - f(s) + f(t) - f(t-s)\big)$.  Thus by the triangle inequality,
$$
\bigg|\E[{\mathcal D}(s) {\mathcal D}(t)] - \E[{\mathcal D}^2(s)]\bigg| \leq \frac{1}{2} \bigg( \big| f(s) \big| + \big| f(t) - f(t-s) \big| \bigg).
$$
(\ref{lipb3}) then follows by combining with the fact that $f$ is Lipschitz with Lipshitz constant at most $\mu^3 \E[S^2]$, while (\ref{lipb4}) follows by combining with (\ref{supf0}) and an additional application of the triangle inequality.  $\Halmos$.
\endproof
\subsubsection{Proof of Lemma\ \ref{boundvarzztop}.}\label{proofzztop}
\proof{Proof:}[Proof of Lemma\ \ref{boundvarzztop}]
The first part of the lemma follows from definitions, the fact that $f$ is Lipschitz with Lipshitz constant at most $\mu^3 \E[S^2]$ by Corollary\ \ref{lipandrest}, and straightforward algebra.  Combining with the basic properties of B.m. and independence of ${\mathcal A}$ and ${\mathcal D}$ then completes the proof of the second part of the lemma.  $\Halmos$.
\endproof
\subsubsection{Proof of Lemma\ \ref{poscovar0}.}\label{proofposcovar0}
\proof{Proof:}[Proof of Lemma\ \ref{poscovar0}]
The lemma follows immediately from Corollary\ \ref{lipandrest}.(\ref{lipb2}), the basic properties of B.m., and independence of ${\mathcal A}$ and ${\mathcal D}$.  $\Halmos$.
\endproof
\subsubsection{Proof of Corollary\ \ref{boundDsupmean}.}\label{proofboundDsupmean}
\proof{Proof:}[Proof of Corollary\ \ref{boundDsupmean}]
It follows from Lemma\ \ref{boundvarzztop} and the stationary increments property of ${\mathcal D}$ that
$\E[\big({\mathcal D}(t) - {\mathcal D}(s)\big)^2] \leq \mu \big(2 c^2_S + 1 \big)(t - s)$ for all $0 \leq s \leq t$.  As for all $0 \leq s \leq t$, 
$$\E\big[\bigg(  \big(2 \mu c^2_S + \mu \big)^{\frac{1}{2}} {\mathcal B}^0(t) -  \big(2 \mu c^2_S + \mu \big)^{\frac{1}{2}} {\mathcal B}^0(s) \bigg)^2\big] =  \mu \big(2 c^2_S + 1 \big) (t - s),$$ we may apply Lemma\ \ref{sudakov} to conclude that 
$$\E[\sup_{0 \leq t \leq T} {\mathcal D}(t)] \leq  \big(2 \mu c^2_S + \mu \big)^{\frac{1}{2}} \E[\sup_{0 \leq t \leq T} {\mathcal B}^0(t)].$$
Combining with Lemma\ \ref{brownprops00}.(\ref{supform}) completes the proof.  $\Halmos$.
\endproof
\subsubsection{Proof of Theorem\ \ref{compareprocess}.}
\proof{Proof:}[Proof of Theorem\ \ref{compareprocess}]
We will prove that the statement of the theorem holds with $f_S = f$, i.e. $f_S(t) = f(t) = V[{\mathcal D}(t)] - \mu c^2_S t$ for all $t \geq 0$, and for any fixed $x \geq \alpha_S$.  That $f$ is continuous and satisfies $\sup_{t \geq 0} |f(t)| \leq \alpha_S$ (from which it follows that $f(t) + 3 x \geq 0$ for all $t \geq 0$) follows immediately from (\ref{supf0}) and Corollary\ \ref{lipandrest}.  Then after applying the definition of $f$, to complete the proof it suffices to demonstrate that for all $0 \leq s < t$,
$$
V[{\mathcal D}'(s), {\mathcal D}'(t)] \leq V[{\mathcal D}(s), {\mathcal D}(t)] + 3 x.$$
Equivalently, after applying definitions and the basic properties of B.m. and the O.U. process, it suffices to demonstrate that for all $0 \leq s < t$,
\begin{equation}\label{compare001}
\mu c^2_S s + \big( f(s) + 3 x \big)^{\frac{1}{2}}\big( f(t) + 3 x \big)^{\frac{1}{2}} \exp\big( - 1.5(t - s) \big) \leq V[{\mathcal D}(s), {\mathcal D}(t)] + 3 x.
\end{equation}
For a fixed $0 \leq s < t$, we treat two cases.  First, suppose $t - s \leq 1$.  In this case, we proceed by first bounding the left-hand-side of (\ref{compare001}) from above, and then bounding the right-hand-side from below.  It follows from (\ref{supf0}), Corollary\ \ref{lipandrest}, and the fact that $x \geq \alpha_S$, that 
\begin{eqnarray*}
\big( f(s) + 3 x \big)^{\frac{1}{2}}\big( f(t) + 3 x \big)^{\frac{1}{2}}
&=& \big(f(s) + 3 x\big) \big( \frac{f(t) + 3 x}{f(s) + 3 x} \big)^{\frac{1}{2}}
\\&=& \big(f(s) + 3 x\big) \big(1 + \frac{f(t) - f(s)}{f(s) + 3 x} \big)^{\frac{1}{2}}
\\&\leq& \big(f(s) + 3 x\big) \big(1 + \frac{x (t-s)}{2 x} \big)^{\frac{1}{2}}.
\\&=& \big(f(s) + 3 x\big) \big(1 + \frac{1}{2}(t-s)\big)^{\frac{1}{2}}.
\end{eqnarray*}
Combining with the fact that $(1 + y)^{\frac{1}{2}} \leq 1 + \frac{1}{2} y$ for all $y \geq 0$, we conclude that 
\begin{equation}\label{compare002}
\big( f(s) + 3 x \big)^{\frac{1}{2}}\big( f(t) + 3 x \big)^{\frac{1}{2}}
\leq
\big(f(s) + 3 x\big) \big(1 + \frac{1}{4}(t-s)\big).
\end{equation}
It follows from a straightforward Taylor-series expansion / error analysis that $\exp\big( - 1.5(t - s) \big) \leq 1 - \frac{3}{4}(t - s)$ for $0 \leq t - s \leq 1$.  Combining with (\ref{compare002}), we conclude that 
\begin{eqnarray*}
&\ &\ \ \big( f(s) + 3 x \big)^{\frac{1}{2}}\big( f(t) + 3 x \big)^{\frac{1}{2}} \exp\big( - 1.5(t - s) \big)
\\&\ &\ \ \ \ \leq\ \ \big( f(s) + 3 x \big) \big(1 + \frac{1}{4}(t-s)\big)\big( 1 - \frac{3}{4}(t - s) \big)
\\&\ &\ \ \ \ =\ \ \big( f(s) + 3 x \big) \big( 1 - \frac{1}{2}(t-s) - \frac{3}{16}(t - s)^2 \big)
\\&\ &\ \ \ \ \leq\ \ \big( f(s) + 3 x \big) \big( 1 - \frac{1}{2}(t-s) \big)
\\&\ &\ \ \ \ =\ \ f(s) + 3 x - \frac{1}{2}\big( f(s) + 3 x\big) (t-s) 
\\&\ &\ \ \ \ \leq\ \ f(s) + 3 x - x (t-s),
\end{eqnarray*}
with the final inequality following from the fact that $f(s) + 3 x \geq 2 x$ by (\ref{supf0}).  Combining with the fact that $\mu c^2_S s + f(s) = \E[{\mathcal D}^2(s)]$, it follows that the left-hand-side of (\ref{compare001}) is at most
\begin{equation}\label{compare003}
\E[{\mathcal D}^2(s)] + 3 x - x (t-s).
\end{equation}
We next bound the right-hand-side of (\ref{compare001}), $V[{\mathcal D}(s), {\mathcal D}(t)] + 3 x$, from below.    
In particular, it follows from Corollary\ \ref{lipandrest} that the right-hand-side of (\ref{compare001}) is at least 
(\ref{compare003}).  This completes the proof of the desired inequality for the case $t - s \leq 1$.
\\\\Next, suppose $t - s > 1$.  We again proceed by bounding each side of (\ref{compare001}).  
It follows from (\ref{supf0}) that 
\begin{equation}\label{compare004}
\big( f(s) + 3 x \big)^{\frac{1}{2}}\big( f(t) + 3 x \big)^{\frac{1}{2}}
\leq
4 x.
\end{equation}
Thus since $1.5(t - s) \geq 1.5$, and $ 4 \exp(-1.5) < 1$, the left-hand-side of (\ref{compare001}) is at most
\begin{equation}\label{compare005}
\mu c^2_S s + 4 x \exp(-1.5) < \mu c^2_S s + x.
\end{equation}
Alternatively, it follows from Corollary\ \ref{lipandrest}, (\ref{supf0}), and the triangle inequality that the right-hand-side of (\ref{compare001}) is at least the right-hand-side of (\ref{compare005}).  This completes the proof of the desired inequality for the case $t - s > 1$, and thus of the theorem as well.  $\Halmos$.
\endproof
\subsubsection{Proof of Corollary\ \ref{dettheoremcor}.}
\proof{Proof:}[Proof of Corollary\ \ref{dettheoremcor}] (\ref{detcor2b}) is proven in \citet{JMM.04} Section 4.2.  We now prove (\ref{detcor2a}).  Let $B' \stackrel{\Delta}{=} B c^{-1}_A$.  Then it follows from Lemma\ \ref{dettheorem1}, and the asymptotic analysis of the supremum of the Gaussian random walk given in \citet{janssen2007cumulants}, in particular Equation 2.18, that 
$$\lim_{B \rightarrow \infty} B^{-2} \log \bigg( \lim_{n \rightarrow \infty} \pr\big( Q^n(\infty) \geq n \big) \bigg)
=  \lim_{B \rightarrow \infty} B^{-2} \log\bigg(1 - \exp\big( - (2 \pi B'^2)^{-\frac{1}{2}} \exp( -\frac{1}{2} B'^2 ) \big) \bigg).
$$
Combining with a straightforward asymptotic analysis completes the proof.
\\\indent Finally, we complete the proof of (\ref{detcor2c}).  Note that for all $x \in \reals$,
$$\pr(S_1 \geq x)\ \ \ \leq\ \ \ \pr\big( \sup_{i \geq 1} S_i \geq x \big)\ \ \ \leq\ \ \ \pr\bigg( \sup_{t \geq 1} \big(c_A {\mathcal B}^0(t) - B t \big) \geq x \bigg).$$ 
The desired claim then follows from the strong Markov property of B.m., and the basic asymptotics of $\Phi$.
$\Halmos$
\endproof

\section*{Acknowledgements.}
The author would like to thank Dimitris Bertsimas, Ton Dieker, David Gamarnik, Kavita Ramanan, and Gennady Samorodnitsky for their helpful discussions and insights.  The author also gratefully acknowledges support from NSF grant no. 1333457.
\bibliographystyle{plainnat}
\bibliography{pwait_bib_7_5_2016}
\end{document}